\documentclass[11pt,a4paper]{article}

\date{}
   
\usepackage{latexsym}
\usepackage{amssymb}

\newtheorem{lemma}{Lemma}
\newtheorem{theorem}{Theorem}

\newcounter{eqclaim}[theorem] 

\def\claim{$$\vtop\bgroup\advance\hsize by -3em\noindent\hspace{-.4em}
\refstepcounter{eqclaim}\ignorespaces\it} \makeatletter
\def\endclaim{\egroup\leqno(\theeqclaim)$$\global\@ignoretrue\rm}
\makeatother

\title{{\bf Transitive orientations in bull-reducible Berge graphs}}

\author{Celina de Figueiredo\thanks{Universidade Federal do Rio de Janeiro.
e-mail: celina@cos.ufrj.br}%
\and%
Fr\'ed\'eric Maffray\thanks{Laboratoire G-SCOP.
e-mail: Frederic.Maffray@g-scop.inpg.fr}%
\and%
Cl\'audia Villela Maciel\thanks{Universidade Federal Fluminense.
e-mail: crvillela@globo.com}%
}

\begin{document}

\maketitle

\begin{abstract}
A bull is a graph with five vertices $r, y, x, z, s$ and five edges
$ry$, $yx$, $yz$, $xz$, $zs$.  A graph $G$ is bull-reducible if no
vertex of $G$ lies in two bulls.  We prove that every bull-reducible
Berge graph $G$ that contains no antihole is weakly chordal, or has a
homogeneous set, or is transitively orientable.  This yields a fast
polynomial time algorithm to color exactly the vertices of such a
graph.
\end{abstract}

%%%%%%%%
%%%%%%%%
\section{Introduction}

A graph is \emph{perfect} if for every induced subgraph $H$ of $G$ the
chromatic number of $H$ is equal to its clique number.  Perfect graphs
were defined by Claude Berge~\cite{ber60}.  The study of perfect
graphs led to several interesting and difficult problems.  The first
one is their characterization.  Berge conjectured that a graph is
perfect if and only if it contains no odd hole and no odd antihole,
where a hole is a chordless cycle of length at least $5$, and an
antihole is the complementary graph of a hole.  It has become
customary to call \emph{Berge graph} any graph that contains no odd
hole and no antihole, and to call the above conjecture the ``Strong
Perfect Graph Conjecture''.  This conjecture was proved by Chudnovsky,
Robertson, Seymour, and Thomas~\cite{CRST}.  A second problem is the
existence of a polynomial-time algorithm to color optimally the
vertices of a perfect graph, solved by Gr\"otschel, Lov\'asz and
Schrijver~\cite{grolovsch84} with an algorithm based on the ellipsoid
method for linear programming.  A third problem is the existence of a
polynomial-time algorithm to decide if a graph is Berge, solved by
Chudnovsky, Cornu\'ejols, Liu, Seymour and
Vu\v{s}kovi\'c~\cite{CCLSV}.  There remains a number of interesting
open problems in the context of perfect graphs, among them the
existence of a combinatorial algorithm to compute the chromatic number
of a perfect graph.

\begin{figure}[htb]
\begin{center}
\begin{picture}(50,50)
    \put(20,0){\circle*{5}}\put(25,-5){$x$}
    \put(0,20){\circle*{5}}\put(-10,18){$y$}
    \put(40,20){\circle*{5}}\put(45,18){$z$}
    \put(0,40){\circle*{5}}\put(-10,38){$r$}
    \put(40,40){\circle*{5}}\put(45,38){$s$}
\put(20,0){\line(-1,1){20}}
\put(20,0){\line(1,1){20}}
\put(0,20){\line(1,0){40}}
\put(0,20){\line(0,1){20}}
\put(40,20){\line(0,1){20}}
\end{picture}
\caption{The bull $r$-$yxz$-$s$.}
\label{fig:bull-labeled}
\end{center}
\end{figure}
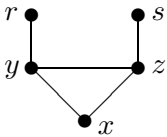

A \emph{bull} is a graph with five vertices $r, y, x, z, s$ and five
edges $ry, yx, yz, xz, zs$; see Figure~\ref{fig:bull-labeled}.  We
will frequently use the notation $r$-$yxz$-$s$ for such a graph.

Bull-free Berge graphs have been much studied as a self-complementary
class of Berge graphs for which first the ``Strong Perfect Graph
Conjecture'' was established by Chv\'atal and Sbihi~\cite{chvsbi87};
subsequently a polynomial-time recognition algorithm for bull-free
Berge graphs was found by Reed and Sbihi~\cite{reesbi}, and further
study of the structure of the class by De Figueiredo, Maffray and
Porto~\cite{defmafpor97,defmafpor01} and Hayward~\cite{hay01} led to a
polynomial-time algorithm to color optimally the vertices of a
bull-free Berge graph~\cite{defmaf04}.

The goal of the present paper is to contribute to the search for a
combinatorial algorithm to compute the chromatic number of a perfect
graph by generalizing the results on the structure of bull-free Berge
graphs~\cite{defmafpor97} to the larger class of bull-reducible Berge
graphs.  A graph $G$ is called \emph{bull-reducible} if every vertex
of $G$ lies in at most one bull of $G$.  Clearly, bull-free graphs are
bull-reducible.  Everett, de~Figueiredo, Klein and Reed~\cite{EFKR}
proved that every bull-reducible Berge graph is perfect.  Although
this result now follows directly from the Strong Perfect Graph Theorem
\cite{CRST}, the proof given in~\cite{EFKR} is much simpler and leads
moreover to a polynomial-time recognition algorithm for bull-reducible
Berge graphs whose complexity is lower than that given for all Berge
graphs in~\cite{CCLSV}.

A graph is called {\it weakly chordal\/} (or ``weakly triangulated'')
if it contains no hole and no antihole.  Hayward {\cite{hay85}} proved
that all weakly triangulated graphs are perfect, and there are very
efficient algorithms to find an exact coloring for weakly chordal
graphs~\cite{hhm89,hss07}.  Given a subset of vertices $S$, a vertex
is said to be {\it partial\/} on $S$ if it has at least one neighbour
and at least one non-neighbour in $S$.  A vertex is {\it impartial on
$S$\/} if it either sees all vertices of $S$ or misses all vertices of
$S$.  A proper subset $H$ of vertices is called {\it homogeneous\/} if
it has at least two vertices and every vertex not in $H$ sees either
all or none of $H$, in other words every vertex not in $H$ is
impartial on $H$.  Notice that if $H$ is a homogeneous set of $G$ then
it is also a homogeneous set of the complement graph $\overline{G}$.
A graph is called \emph{transitively orientable} if it admits a
\emph{transitive orientation\/}, i.e., an orientation of its edges
with no circuit and with no $P_3$ $abc$ with the orientation
$\vec{ab}$ and $\vec{bc}$.  Here we prove:

\begin{theorem}\label{main1}
Let $G$ be a bull-reducible Berge graph that contains no antihole.  
Then $G$ either is weakly chordal, or has a homogeneous set, or is
transitively orientable.
%transitively orientable, or contains a homogeneous set.
\end{theorem}
Using Theorem~\ref{main1}, we can devise a polynomial-time algorithm
that colors the vertices of any bull-reducible Berge graph that
contains no antihole.  This question will be addressed in
Section~\ref{s:co}. 

%The rest of this manuscript consists in a 
%sketch of the proof of this theorem. 
%

%\begin{theorem}\label{main2}
%Let $G$ be a graph in ${\cal B}$.  If $G$ contains a hole and an
%antihole then $G$ has a homogeneous pair.  Moreover, the two graphs
%obtained by the decomposition of $G$ along this homogeneous pair are
%in ${\cal B}$\footnote{THE PROOF OF THIS IS STILL MISSING.}.
%\end{theorem}
%%
%\emph{Proof.} {\it For the first part, try to imitate Chv\'atal and 
%Sbihi. For the second part, try to imitate ``Optimizing bull-free 
%perfect graphs''.  \ldots}
\

A \emph{wheel} is a graph made of an even hole of length at least $6$
plus a vertex that sees all vertices of this hole.  A \emph{double
broom} is a graph made of a $P_4$
%(called the central $P_4$ of the double broom), 
plus two non-adjacent vertices $a, a'$ that see all vertices of the
$P_4$, plus a vertex $b$ that sees only $a$ and a vertex $b'$ that
sees only $a'$.  A \emph{lock} is a graph with six vertices such that
the first four induce a hole, the fifth one is adjacent to the first
four, and the sixth one is adjacent to two adjacent vertices of the
hole only.
%A \emph{domino} is a graph with six vertices obtained from a $C_6$ by
%adding one edge between two vertices at distance $3$ along the $C_6$.
Given a graph $F$, a \emph{spiked $F$} is a graph that consists in a
copy of $F$ plus two additional vertices $a,b$ such that $b$ has no
neighbour in $F$ and $a$ is adjacent to every vertex of $V(F)\cup
\{b\}$.  Let us use the notation $F_1$, $F_2$ for the following two
types of graphs: $F_1$ stands for the bull, and $F_2$ for the lock.
Let ${\cal B}$ be the class of bull-reducible Berge graphs that
contain no wheel, no double broom, and no spiked $F_j$ ($j=1,2$).

\begin{theorem}\label{to}
Let $G$ be a graph in ${\cal B}$.  If $G$ contains a hole of length at
least six and no antihole then $G$ is transitively orientable.
\end{theorem}

\begin{lemma}[\cite{EFKR,FMV}]
\label{lem:homoset}
Let $G$ be a bull-reducible $C_5$-free graph.  If $G$ contains a wheel
or a double broom
then $G$ has a homogeneous set. 
\end{lemma}

%\begin{lemma}[\cite{FMV}]
%Let $G$ be a bull-reducible $C_5$-free graph.  If $G$ contains a
%double broom then $G$ has a homogeneous set.  
%\end{lemma} 

\begin{lemma}\label{lem:spkb} 
Let $G$ be a bull-reducible $C_5$-free graph.  If $G$ contains a
spiked $F_j$ for any $j=1, 2$, then $G$ has a homogeneous set.
\end{lemma}

Lemma~\ref{lem:spkb} is proved in Section~\ref{sec:lemmas}, and
Theorem~\ref{to} is proved in Section~\ref{s:to}.  We can see
immediately how to obtain a proof our main Theorem.

\

{\it Proof of Theorem~\ref{main1}.} Let $G$ be a bull-reducible Berge
graph containing no antihole.  If $G$ contains a wheel, a double
broom, or a spiked $F_i$ ($i=1,2$) then $G$ has a homogeneous set by
Lemmas~\ref{lem:homoset} and~\ref{lem:spkb}.  So we may assume that
$G$ is a graph in the class ${\cal B}$.  If $G$ contains no hole of
length at least six then it is weakly chordal.  So we may assume that
$G$ contains a hole of length at least six.  Then Theorem~\ref{to}
implies that $G$ is transitively orientable.  $\Box$

\section{Some lemmas}
\label{sec:lemmas}

\emph{Proof of Lemma~\ref{lem:spkb}.} Suppose that $G$ has an induced
subgraph $S$ that is a spiked $F_j$ for some $j=1, 2$.  Let $W$ be the
set of vertices that induces the $F_j$ contained in $S$; let $b$ be
the vertex of $S$ that misses every vertex of $W$; and let $a$ be the
vertex of $S$ that sees all of $W\cup\{b\}$.  Let $W$ have vertices
$u_1, \ldots, u_{|W|}$ with the following notation.  If $W$ induces a
bull ($F_1$), then it is $u_1$-$u_2u_5u_3$-$u_4$.  If $W$ induces a
lock ($F_2$), it has edges $u_1u_2, u_2u_3, u_3u_4, u_4u_1$, $u_5u_i$
($i=1, \ldots, 4$), $u_1u_6, u_2u_6$.  In either case, we define
additional sets of vertices as follows.  Let $T$ be the set of
vertices of $G-W$ that see all vertices of $W$.  Let $Z$ be the set of
vertices of $G-W$ that see none of $W$.  Let $P$ be the set of
vertices of $G-W$ that have a neighbour and a non-neighbour in $W$.
Clearly, $W, T, Z, P$ form a partition of $V(G)$.

For any $p\in P$, say that three vertices $u,v,w\in W$ form a
\emph{blue triple} if they induce a subgraph with only one edge, say
the edge $uv$, and $p$ sees $u$ and misses $v,w$.  Note that if there
is such a triple and $p$ misses any $t\in T$ then $p$-$uvt$-$w$ is a
bull in $G$, and we call any such bull a ``blue bull''.  Thus,
\begin{claim}\label{2blue}
If $p$ has two blue triples, then it sees all of $T$.
\end{claim}

For any $p\in P$, say that a chordless path $u$-$v$-$w$ of three
vertices of $W$ is \emph{red} if $p$ sees $u,v$ and misses $w$.  Note
that if there is such a path and $p$ sees any $z\in Z$ then
$z$-$puv$-$w$ is a bull in $G$, and we call any such bull a ``red
bull''.  Thus,
\begin{claim}\label{2red}
If $p$ has two red paths, then it misses all of $Z$.
\end{claim}

We will need a classification of the vertices of $P$.  Let $p$ be a
vertex in $P$; then, using only the fact that $W\cup\{p\}$ induces a
subgraph that contains no $C_5$ and at most one bull, it is a routine
matter to establish that $p$ must be of exactly one of the types
listed in the following table:
\small
\begin{center}
\begin{tabular}{|c|c|c|c|}
    \hline\hline
$W$  & $N(p)\cap W$ (up to symmetry) & number of & number of \\
 & & blue triples & red paths \\ \hline\hline
$F_1$ & $\{u_1\}$ or $\{u_1, u_3\}$ & $1$ & $0$ \\ \cline{2-4}
 & $\{u_1, u_2, u_5\}$ or $\{u_1, u_2, u_4, u_5\}$ & $0$ & $1$ \\ \cline{2-4}
  &  $\{u_2, u_5\}$ & $2$ & $1$ \\ \cline{2-4}
  &  $\{u_1, u_2, u_3\}$ & $1$ & $2$    \\ \hline\hline
  $F_2$ & $\{u_1\}$ & $1$ & $0$ \\ \cline{2-4} 
  & $\{u_3\}$ or $\{u_1, u_3\}$ & $2$ & $0$ \\ \cline{2-4} 
  & $\{u_3, u_5\}$, $\{u_1, u_3, u_5\}$ or $\{u_1, u_3, u_5, u_6\}$ &
  $\ge 1$ & $1$ \\ \cline{2-4}
  %%%%%%%%%%%%%%%%%%%%%%%%%%%%%%%%
  & other possible types  & -- & $\ge 2$  \\   \hline\hline
\end{tabular}
\end{center}

The ``other possible types'' for $F_2$ are (up to symmetry): \\ 
$\{u_1, u_2\}$, $\{u_1, u_6\}$, $\{u_1, u_2, u_4\}$, $\{u_1, u_2,
u_6\}$, $\{u_1, u_3, u_4\}$, $\{u_1, u_3, u_6\}$, $\{u_1, u_5, u_6\}$,
$\{u_3, u_4, u_5\}$, $\{u_1, u_2, u_3, u_4\}$, $\{u_1, u_2, u_4,
u_6\}$, $\{u_1, u_2, u_5, u_6\}$ $\{u_1, u_3, u_4, u_5\}$, $\{u_3,
u_4,$ $u_5,$ $u_6\}$, $\{u_1, u_2, u_3, u_4, u_5\}$, $\{u_1, u_2, u_3,
u_4, u_6\}$, $\{u_1, u_2, u_4, u_5, u_6\}$, and $\{u_1, u_3, u_4,$
$u_5,$ $u_6\}$.

\normalsize%

Let $P_1$ be the set of vertices $p$ of $P$ such that $N(p)\cap W$ is
a stable set, and let $P_2= P-P_1$.  We claim that:
\begin{claim}\label{tp1}
Every vertex of $T$ sees every vertex of $P_1$.
\end{claim}
%
% \section*{Appendix}
% \subsection*{Proof of claim \ref{tp1} in Lemma \ref{lem:spkb}}
For suppose on the contrary that there are non-adjacent vertices $t\in
T$ and $p\in P_1$.  By the remark, $p$ has no red path.  By
(\ref{2blue}), $p$ has at most one blue triple.  The table shows two
types that satisfy these conditions, on the 1st line of $F_1$ and the
1st line of $F_2$.  In either case there is a blue bull.  If $W$
induces an $F_1$ then $W$ and the blue bull intersect, a
contradiction.  If $W$ induces an $F_2$ then $p$-$u_1u_6u_2$-$u_3$ is
a second bull containing $p$, a contradiction.  Therefore (\ref{tp1})
holds.\\

Let $A$ be the set of those vertices of $T$ that have a neighbour in
$Z$.  Clearly $b\in Z$ and $a\in A$.  We claim that:
\begin{claim}\label{ap}
Every vertex of $A$ sees every vertex of $P$.
\end{claim}
%
% \subsection*{Proof of claim \ref{ap} in Lemma \ref{lem:spkb}}
For suppose on the contrary that some vertex $t\in A$ misses some
vertex $p\in P$.  Up to renaming vertices we may assume that $t=a$.
Since $A\subseteq T$, and by (\ref{tp1}), we have $p\in P_2$, so there
is at least one red path for $p$.  Suppose that $p$ sees $b$.  Then
(\ref{2red}) implies that there is exactly one red path, and there is
a red bull.  Thus we must have zero blue triple for $p$.  The table
shows one type that satisfies these conditions, on the 2nd line of
$F_1$.  In this case, $W$ and the red bull are two intersecting bulls,
a contradiction.  Thus $p$ misses $b$.  Say that an edge $uv$ with
$u,v\in W$ is a \emph{switch} if $p$ sees $u$ and misses $v$.  Note
that if there is such an edge then $b$-$avu$-$p$ is a bull in $G$ (a
``switch bull'').  Thus, there must be at most one switch.  In fact
there is a switch since $W$ is connected and both $W\cap N(p)$ and
$W-N(p)$ are non empty; so there is a switch bull, and consequently
there must be zero blue triple for $p$.  If $W$ induces an $F_2$, then
$W$ is $2$-connected and so there are two switches, a contradiction.
If $W$ induces an $F_1$ then the switch bull and $W$ itself are two
intersecting bulls, a contradiction.  Therefore (\ref{ap}) holds.\\

Let $X$ be the set of vertices $x$ of $Z$ such that there exists in
$G$ a path $x_0$-$x_1$-$\cdots$-$x_k$ with $x_0\in P$, $k\ge 1$, $x_1,
\ldots, x_k\in Z$, and $x=x_k$.  We claim that:
\begin{claim}\label{ax}
Every vertex of $A$ sees every vertex of $X$.
\end{claim}
%
% \subsection*{Proof of claim \ref{ax} in Lemma \ref{lem:spkb}}
For suppose that some vertex $t\in A$ misses some vertex $x\in X$.  Up
to renaming vertices we may assume that $t=a$.  By the definition of
$X$ there is a path $x_0$-$\cdots$-$x_k$ with $x_0\in P$, $x_1,
\ldots, x_k\in Z$, and $x=x_k$.  We may assume that $k$ is minimal, so
this path is chordless.  By (\ref{ap}), $a$ sees $x_0$.  Since $W$ is
connected in $\overline{G}$, there are non-adjacent vertices $w,w'\in
W$ such that $x_0$ sees $w$ and misses $w'$.  Suppose that $k=1$.
Then $x_1$-$x_0wa$-$w'$ is a bull in $G$, so there must be only one
such pair $w,w'$.  %
%This is possible only if $W$ is not $2$-connected
% in $\overline{G}$, which means that $W$ induces an $F_1$ or $F_2$.
When $W$ induces an $F_1$, $W$ itself is a second bull containing $w$,
a contradiction.  So let $W$ induce $F_2$.  If $x_0\in P_1$ then the
table shows that $N(x_0)\cap W\subseteq \{u_1, u_3\}$, and either
$x_0$-$u_1u_6u_2$-$u_3$ (if $x_0$ misses $u_3$) or
$x_0$-$u_3u_5u_2$-$u_6$ (if $x_0$ sees $u_3$) is a second bull
containing $x_0$, a contradiction.  If $x_0\in P_2$, then there is a
red path in $W$, so there is a red bull with $x_0, x_1$, which is a
second bull containing $x_0$, a contradiction.  Now suppose that $k\ge
2$.  The minimality of $k$ implies that $a$ sees $x_{k-2}$ and
$x_{k-1}$.  Let $w'$ be any vertex in $W-N(x_{k-2})$.  Then
$w'$-$ax_{k-2}x_{k-1}$-$x_k$ is a bull containing $x_0$ for each
choice of $w'$, which is a contradiction if $|W-N(x_{k-2})|\ge 2$.  So
we must have $|W-N(x_{k-2})|=1$, which implies $k=2$ and,
by~(\ref{tp1}), $x_0\in P_2$; but then there is a red path and a red
bull with $x_0, x_1$, which is a second bull containing $x_0$, a
contradiction.  Therefore (\ref{ax}) holds.\\

Let $Y$ be the set of vertices of $T-A$ such that there exists in
$\overline{G}$ a path $x$-$y_1$-$\cdots$-$y_\ell$ with $x\in P\cup X$,
$\ell\ge 1$, $y_1, \ldots, y_\ell\in T-A$ and $y=y_\ell$.  We claim that:
\begin{claim}\label{ay}
Every vertex of $A$ sees every vertex of $Y$.
\end{claim}
%
% \subsection*{Proof of claim \ref{ay} in Lemma \ref{lem:spkb}}
For suppose that some vertex $t\in A$ misses a vertex $y\in Y$.  Up to
renaming vertices we may assume that $t=a$.  By the definition of $Y$
and $X$, there is a sequence of vertices $x_0, x_1, \ldots, x_k$,
$y_1, \ldots, y_\ell$ such that $x_0$-$x_1$-$\cdots$-$x_k$ is a path
in $G$ and $x_k$-$y_1$-$\cdots$-$y_\ell$ is a path in $\overline{G}$
such that $x_0\in P$, if $k\ge 1$ then $x_1, \ldots, x_k\in X$,
$\ell\ge 1$, $y_1, \ldots, y_\ell\in Y$ and $y=y_\ell$.  We may assume
that this sequence is minimal with these properties.  Note that $y_1,
\ldots, y_\ell$ miss $b$ since they are in $Y$.  If $\ell\ge 2$, then,
by the minimality of the sequence, $a$ sees $y_{\ell-1}$, and so
$b$-$ay_{\ell-1}w$-$y_\ell$ is a bull for each $w\in W$, a
contradiction.  So $\ell=1$.  If $k\ge 2$, then, by the minimality of
the sequence, $y_1$ sees $x_1$, which contradicts the definition of
$Y$.  So $k\le 1$.  \\
Suppose that $k=0$.  Then $y_1$ misses $x_0$, which implies, by
(\ref{tp1}), that $x_0\in P_2$ (so there is a red path) and, by
(\ref{2blue}), that $x_0$ has at most one blue triple.  If $x_0$
misses $b$, then $b$-$ax_0w$-$y_1$ is a bull for each neighbour $w$ of
$x_0$ in $W$, which is possible only if there is only one such $w$, so
$x_0\in P_1$, a contradiction.  Thus $x_0$ sees $b$.  Then
(\ref{2red}) implies that there is at most one, and so exactly one,
red path for $x_0$, and there is one red bull with $x_0, b$.
Consequently there must be zero blue triple.  The table shows one type
that satisfies these conditions, on the 2nd line $F_1$.  In this case,
$W$ and the red bull are two intersecting bulls, a contradiction.  \\
Suppose that $k=1$.  By the minimality, $y_1$ sees $x_0$.  Since $x_0$
has a neighbour $x_1$ in $Z$, (\ref{2red}) implies that there is at
most one red path for $x_0$.  If there is one red path $u$-$v$-$w$,
then there is a red bull $x_1$-$x_0uv$-$w$ and a second bull
$x_1$-$x_0uy_1$-$w$, a contradiction.  So there is no red path.  The
table shows three types that satisfy this condition: on the 1st line
of $F_1$, and the 1st and 2nd lines of $F_2$.  In either case let $w,
w'$ be any two non-adjacent vertices of $W$ such that $x_0$ sees $w$
and misses $w'$.  Such a pair exists since $W$ is connected in
$\overline{G}$.  Then $x_1$-$x_0wy_1$-$w'$ is one bull.  So there must
be only one such pair $w,w'$.  This implies that $W$ must not be
$2$-connected in $\overline{G}$.  If $W$ is $F_1$, then $N(x_0)\cap W$
is one of $\{u_1\}$ or $\{u_1, u_3\}$ (1st line) and there are also
two pairs ($u_1, u_4$ and $u_1, u_5$) like $w,w'$.  So it must be that
$W$ induces an $F_2$ and $N(x_0)\cap W$ is one of $\{u_1\}$,
$\{u_3\}$, or $\{u_1, u_3\}$.  If it is $\{u_3\}$ then again there are
two pairs ($u_3, u_1$ and $u_3, u_6$) like $w,w'$; if it is $\{u_1\}$,
then $x_0$-$u_1u_6u_2$-$u_3$ is a second bull containing $x_0$, and if
it is $\{u_1, u_3\}$ then $x_0$-$u_3u_5u_2$-$u_6$ is a second bull, a
contradiction.  Therefore (\ref{ay}) holds.  \\

Now $V(G)$ can be partitioned into the set $H = W\cup P\cup X\cup Y$
and the sets $T-Y$ and $Z-X$.  We claim that:
\begin{claim}\label{ht}
Every vertex $h\in H$ sees every vertex $t\in T-Y$.
\end{claim}
If $h\in W$, this is by the definition of $T$.  If $h\in P\cup X\cup
Y$ and $t\in A$ this is by~(\ref{ap}), (\ref{ax}) and~(\ref{ay}).   
If $h\in P\cup X\cup Y$ and $t\in T-Y-A$ this is by the definition of
$Y$.  Thus (\ref{ht}) holds.

Next we claim that:
\begin{claim}\label{hz}
Every vertex $h\in H$ misses every vertex in $Z-X$.
\end{claim}
If $h\in W$ this is by the definition of $Z$.  If $h\in P\cup X$ this
is by the definition of $X$.  If $h\in Y$ this is by the definition of
$Y$ ($\subseteq T-A$). Thus (\ref{hz}) holds.

Now it follows from (\ref{ht}), (\ref{hz}) and the fact that $T-Y$ is
not empty (since it contains $a$) that $H$ is a homogeneous set of
$G$, which concludes the proof of the lemma.  $\Box$

\
 
We finish this section by recalling a useful lemma. 

\begin{lemma}[\cite{EFKR}]
\label{lem:neighbourhood}
Let $G$ be a bull-reducible odd hole-free graph, and let $C$ be a
shortest even hole of length at least six in $G$, with its vertices
labelled alternately ``odd'' and ``even''.  Let $v$ be any vertex in
$V(G)\setminus V(C)$.  Then $v$ satisfies exactly one of the following
conditions:
\begin{itemize}
\item
$N(v)\cap V(C) = V(C)$, so $C$ and $v$ form a wheel;
\item
$N(v)\cap V(C)$ consists in either all even vertices and no odd vertex
of $C$ or all odd vertices and no even vertex of $C$;
\item
$N(v)\cap V(C)$ consists in either one, or two consecutive or three
consecutive vertices of $C$;
\item
$N(v)\cap V(C)$ consists in two vertices at distance $2$ along $C$;
\item
$C$ has length $6$ and $N(v)\cap V(C)$ consists in four vertices such
that exactly three of them are consecutive.
\end{itemize}
\end{lemma}

%%%%%%%%
%%%%%%%%
\section{Transitive orientations}
\label{s:to}

In a graph $G$, for any set $B\subseteq V(G)$, let $P(B)$ be the set
of vertices of $V(G)-B$ that have a neighbour and a non-neighbour in
$B$.  Let us say that a graph $G$ has a box partition if its vertex
set can be partitioned into non-empty subsets, called boxes, with the
following properties:
\it
\begin{enumerate}
\item[(i)]
Every box is labelled either odd or even, and there is no edge between
two boxes that have the same label.
\item[(ii)]
Every box induces a connected subgraph of $G$.
\item[(iii)]
Boxes are either central or peripheral, and there are at least six
central boxes.
\item[(iv)]
For each central box $B$ there are four auxiliary vertices $a_B, b_B,
c_B, d_B \in V(G)-B$, such that $a_B$ and $c_B$ see all of $B$ and
miss all of of $P(B)$, $b_B$ and $d_B$ miss all of $B$, and the four
auxiliary vertices induce a graph with two edges $a_Bb_B,
c_Bd_B$.
\item[(v)]
For each peripheral box $B$ there are two auxiliary vertices $a_B,
b_B, \in V(G)-B$, such that $a_B$ sees all of $B$ and misses all of of
$P(B)$, $b_B$ misses all of $B$, and the two auxiliary vertices are
adjacent.
\item[(vi)]
For any two adjacent vertices $u,v$ in a peripheral box $B$ the sets
$N(u)\cap P(B)$ and $N(v)\cap P(B)$ are comparable by inclusion.
\item[(vii)]
A box $B$ does not contain a chordless path $u$-$v$-$w$ such that
there are adjacent vertices $x,y\in P(B)$ such that $x$ sees $u$ and
misses $v,w$ and $y$ sees $u,v$ and misses $w$.
\end{enumerate}
\rm

% LEMMA "Moreitems"
\begin{lemma}\label{lem:moreitems}
%%Let $G$ be any bull-reducible, $C_5$-free graph, and let $B\subseteq
%%V(G)$ be a set that satisfies~(iv) or~(v).  Suppose that there are
%%vertices $u,v\in B$ and $x,y\in P(B)$ such that $x$ ses $u$ and misses
%%$v$ and $y$ sees $v$ and misses $u$.  Then $u,v$ are adjacent and $G$
%%contains a $\overline{C}_6$.
Let $G$ be any bull-reducible, $C_5$-free graph, and let $B\subseteq
V(G)$ be a set that satisfies~(iv) or~(v).  Suppose that there are
vertices $u,v\in B$ and $x,y\in P(B)$ such that $x$ sees $u$ and misses
$v$ and $y$ sees $v$ and misses $u$. 
Then there is no $P_4$ induced by $x$, $u$, $v$, $y$.
If $u$, $v$ are adjacent, then $G$ contains a $\overline{C}_6$.
% 3.  If a set $B\subseteq V(G)$ satisfies~(iv), there is no $P_3$
% $u$-$v$-$w$ in $B$ with a vertex $x\in P(B)$ that sees $u,v$ and
% misses $w$. \\
\end{lemma}
%
% \subsection*{Proof of Lemma \ref{lem:moreitems}}
 \emph{Proof.}
Suppose $x$, $u$, $v$, $y$ induce a $P_4$.
Since $B$ satisfies~(iv) or~(v), then, with the same
notation, we know that $a_B$ sees $u, v$ and misses $x,y$, and $b_B$
misses $u, v$.  
% Suppose that $x,y$ are not adjacent.  Then
$x$-$ua_Bv$-$y$ is a bull.  If $b_B$ misses $x$ then $b_B$-$a_Bvu$-$x$
is a second bull, which intersects the first.  So $b_B$ sees $x$, and
similarly it sees $y$.  But then $b_B, x, u, v, y$ induce a $C_5$, a
contradiction.
% So $x,y$ are adjacent.  
Now suppose $u$, $v$ are adjacent. This implies $x,y$ are adjacent.
Then $b_B$ sees at least one
of $x,y$, for otherwise $b_B$-$a_Buv$-$y$ and $b_B$-$a_Bvu$-$x$ are
two intersecting bulls.  Then $b_B$ sees both $x,y$, for if it sees
only one, say $x$, then $a_B, b_B, x, y, v$ induce a $C_5$.  Thus
$a_B, b_B, u, v, x, y$ induce a $\overline{C}_6$.  This completes the
proof of the Lemma.  $\Box$

% A CONSERVER, UTILE
%
% If a set $B\subseteq V(G)$ satisfies~(\ref{it:auxi}) and $G$ contains
% no $\overline{C}_6$, then for any two adjacent vertices $u,v \in B$
% the sets $N(u)\cap P(B)$ and $N(v)\cap P(B)$ are comparable by
% inclusion.
% 
% %
% Now suppose that~(\ref{it:noc6}) fails: There are vertices $x,y\in
% P(B)$ such that $x$ sees $u$ and misses $v$ and $y$ sees $v$ and
% misses $u$.  By~(\ref{it:nop4}), $x$ sees $y$.  Then $b_B$ sees one of
% $x,y$, for otherwise $b_B$-$a_Buv$-$y$ and $b_B$-$a_Bvu$-$x$ are tow
% intersecting bulls.  Then $b_B$ sees both $x,y$, for if it sees only
% one, say $x$, then $b_B, x, y, v, a_B$ induce a $C_5$.  Now $a_B, b_B,
% u, v, x, y$ induce a $\overline{C}_6$.   

\

Let us say that a vertex $x$ in a graph $G$ is \emph{sensitive} if
there exist six vertices $u_1, \ldots, u_6$ of $G-x$ with edges
$u_iu_{i+1}$ ($i=1, \ldots, 5$) and possibly $u_1u_6$ (so that they
induce a $P_6$ or $C_6$) such that $x$ is adjacent to $u_2$ and $u_3$
and not to $u_1, u_4, u_5, u_6$, and $G-x$ contains a hole.   

% A \emph{long bull} is a graph with seven vertices $x, u_1, \ldots,
% u_6$ and edges $u_iu_{i+1}$ ($i=1, \ldots, 5$), $xu_2, xu_3$, and
% possibly the edge $u_1u_6$ (so the $u_i$'s induce a $P_6$ or $C_6$).
% Vertex $x$ is called the \emph{nose} of the long bull.

%Say that a bull $a$-$bxc$-$d$ is \emph{cyclic} if the path
%$a$-$b$-$c$-$d$ extends to a hole.  Note that in that case this hole
%is even and has length at least $6$, and Lemma~\ref{lem:neighbourhood}
%%the Neighbourhood Lemma
%implies that $x$ sees only $b$ and $c$ along that hole.  Say that a
%bull $a$-$bxc$-$d$ is \emph{pericyclic} if there is a hole $C$ of
%length $\ge 6$ and a chordless path $P=p_1$-$\cdots$-$p_k$ such that:
%$p_1$ has a neighbour $q$ in $C$, there is no edge between $P\setminus
%\{p_1, p_2\}$ and $C$, $P+q$ contain $a$-$b$-$c$-$d$, vertex $x$ has
%no neighbour in $C\cup P$ other than $b,c$, and $C\cup P$ contains no
%triangle.

%%%%%%%%
\begin{theorem}\label{nocycb}
Let $G$ be a graph in ${\cal B}$.  If $G$ contains a hole of length at
least $6$ and $G$ has no sensitive vertex, 
%cyclic bull and no pericyclic bull, 
then $G$ admits a box partition.
\end{theorem}
\emph{Proof.} Let $\ell$ be the length of a shortest even hole of
length at least $6$ in $G$.  So there exist $\ell$ non-empty disjoint
subsets $V_1, \ldots, V_\ell$ such that each vertex in $V_i$ sees
every vertex in $V_{i-1}\cup V_{i+1}$ and misses every vertex in
$V_{i+2}\cup V_{i+3}\cup\cdots\cup V_{i-2}$ (with subscripts modulo
$\ell$).  Put $V^*_1=V_1\cup V_3 \cup\cdots\cup V_{\ell-1}$,
$V^*_2=V_2\cup V_4 \cup\cdots\cup V_\ell$, and $V^*=V^*_1\cup V^*_2$.
We may assume that $V^*$ is maximal with this property.  Let us then
define the following subsets of vertices:
\begin{itemize}
\item
Let $A^*_1$ be the set of vertices that see all of $V^*_2$ and miss
all of $V^*_1$;
\item
Let $A^*_2$ be the set of vertices that see all of $V^*_1$ and miss
all of $V^*_2$;
\item
For $i=1, \ldots, \ell$, let $X_i$ be the set of all vertices not in
$A^*_1\cup A^*_2$ that see all of $V_{i-1}\cup V_{i+1}$ and miss all
of $V_{i-2}\cup V_{i+2}$;
\item
$D_i = V_i\cup X_i$;
\item
$D^*_1= D_1\cup D_3\cup\cdots\cup D_{\ell-1}$, $D^*_2= D_2\cup
D_4\cup\cdots\cup D_\ell$;
\item
$C^*_1 = D^*_1\cup A^*_1$, $C^*_2 = D^*_2\cup A^*_2$;
\item
$Z = V(G) - (D^*_1\cup D^*_2\cup A^*_1\cup A^*_2)$.
\end{itemize}
Clearly, the sets $D_1,\ldots, D_\ell, A^*_1, A^*_2, Z$ form a
partition of $V(G)$.  Note that subscripts on the starred sets are
modulo $2$, while subscripts on the unstarred sets are modulo $\ell$.
From now on we reserve the letter $v_i$ for an arbitrary vertex in
$V_i$ ($i=1,\dots, \ell$).  Let us establish a number of useful facts.
\begin{claim} 
If any $X_i$ is non-empty, then $\ell=6$.  Every vertex of $X_i$ has a
neighbour in $V_{i+3}$.  If a vertex of $X_i$ sees all of $V_{i+3}$
then it has a neighbour in~$V_i$.
\label{xivi}
\end{claim}
%
% \subsection*{Proof of claim \ref{xivi} in Theorem \ref{nocycb}}
For simpler notation put $i=3$.  Let $x$ be any vertex of $X_3$.  So
$x$ sees all of $V_2\cup V_4$ and misses all of $V_1\cup V_5$.  Then
$x$ must have a neighbour in $V_6\cup \cdots\cup V_\ell$, for
otherwise we could add $x$ to $V_3$, which would contradict the
maximality of $V^*$.  Let $h$ be the smallest index such that $x$ has
a neighbour $y$ in $V_h$ with $6\le h\le \ell$.  If $h\ge 7$, then
$\{x, v_4, \ldots, v_{h-1}, y\}$ induces a hole of length $h-2$, with
$5\le h-2\le l-2$, which contradicts $G$ being Berge (if $h$ is odd)
or the definition of $\ell$ (if $h$ is even).  So $h=6$.  Suppose
$\ell\ge 8$.  Then we can apply Lemma~\ref{lem:neighbourhood} to the
hole induced by $\{v_1, v_2, v_3, v_4, v_5, y, \ldots, v_\ell\}$ and
to $x$, which implies that $x$ sees every $v_j$ with even $j\neq 6$
and misses every $v_j$ with odd $j$.  Then applying
Lemma~\ref{lem:neighbourhood} to the hole induced by $\{v_1, \ldots,
v_\ell\}$ implies that $x$ also sees every $v_6\in V_6$.  But then we
have $x\in A^*_1$, which contradicts the definition of $X_3$.  Thus
$\ell=6$.  Now if $x$ also sees all of $V_6$ and none of $V_3$, then
$x$ must be in $A_1$, which contradicts the definition of $X_3$.  So
if $x$ sees allf of $V_6$ it has a neighbour in $V_3$.  Therefore
(\ref{xivi}) holds.
\begin{claim}\label{nodidj}
For $i,j$ of the same parity, there is no edge between $D_i$ and
$D_j$.
\end{claim}
%
% \subsection*{Proof of claim \ref{nodidj} in Theorem \ref{nocycb}}
For if $\ell\ge 8$, this follows immediately from the fact that
$D_i=V_i$ and $D_j=V_j$.  Now let $\ell=6$ and suppose up to symmetry
that there is an edge $xy$ with $x\in D_1$ and $y\in D_3$.  Since $x$
has a neighbour in $D_3$ we have $x\notin V_1$, so $x\in X_1$; and
then, by (\ref{xivi}), we have $\ell=6$ and $x$ has a neighbour
$u_4\in V_4$.  Likewise, $y$ is in $X_3$ and has a neighbour $u_6\in
V_6$.  If $x$ has a non-neighbour $w_4\in V_4$ and $y$ has a
non-neighbour $w_6\in V_6$ then $\{x, y, w_4, v_5, w_6\}$ induces a
$C_5$, a contradiction.  So we may assume, up to symmetry, that $x$
sees all of $V_4$.  Then (\ref{xivi}) implies that $x$ has a neighbour
$w_1\in V_1$.  So we find a bull $w_1$-$xyu_4$-$v_5$.  If $y$ has a
neighbour $w_3\in V_3$, then we find a second bull $w_3$-$yxu_6$-$v_5$
containing $x$, a contradiction.  So $y$ has no neighbour in $V_3$,
and, by~(\ref{xivi}), $y$ has a non-neighbour $w_6\in V_6$.  But then
we find a second bull $v_5$-$w_6w_1x$-$y$, a contradiction.
Therefore~(\ref{nodidj}) holds.
\begin{claim}\label{noaidi}
There is no edge between $A^*_1$ and $D^*_1$.\\
There is no edge between $A^*_2$ and $D^*_2$.
\end{claim}
%
% \subsection*{Proof of claim \ref{noaidi} in Theorem \ref{nocycb}}
For suppose, up to symmetry, that a vertex $a$ in $A^*_1$ sees a
vertex $x_1$ in $D_1$.  By the definition of $A^*_1$, $x_1$ is in
$X_1$, so $\ell=6$ by (\ref{xivi}), and $x_1$ has a neighbour $w_4\in
V_4$.  If $x_1$ sees any $v_1\in V_1$ then $v_1$-$x_1aw_4$-$v_3$ and
$v_1$-$x_1aw_4$-$v_5$ are two intersecting bulls, a contradiction.  If
$x_1$ misses every $v_1\in V_1$, then by (\ref{xivi}) it misses some
$u_4\in V_4$.  But then $v_1$-$v_2x_1a$-$u_4$ and
$v_1$-$v_6x_1a$-$u_4$ are two intersecting bulls, a contradiction.
Therefore (\ref{noaidi}) holds.
%
%
% APPARENTLY THE FOLLOWING CLAIM IS NO LONGER USEFUL
%
% \begin{claim}\label{diimpa}
% Every vertex of $D_i$ is impartial on each component of $A^*_1$ and of
% $A^*_2$.
% \end{claim}
% 
% For suppose that for some $x\in D_i$ there are adjacent vertices $a,b
% \in A^*_1$ such that $x$ sees $a$ and misses $b$.  By (\ref{noaidi})
% and up to symmetry we may assume $i=2$.  By the definition of $D_2$ we
% may also assume that $\ell=6$ and that $x\in X_2$.  If $x$ sees every
% $v_5 \in V_5$ then, by (\ref{xivi}), it has a neighbour $w_2\in V_2$;
% but then $v_5, x, w_2, b, v_4$ induce a $C_5$.  So $x$ has a
% non-neighbour $w_5\in V_5$.  Then $x$-$abv_4$-$w_5$ and
% $x$-$abv_6$-$w_5$ are two intersecting bulls, a contradiction.
% Therefore (\ref{diimpa}) holds.
%
%
\begin{claim}\label{zv12}
Each $z\in Z$ misses all of $V^*_1$ or all of $V^*_2$, and there is at
most one $i\in\{1, \ldots, l\}$ such that $z$ sees all of $V_i$.
\end{claim}
%
% \subsection*{Proof of claim \ref{zv12} in Theorem \ref{nocycb}}
To prove the first part of the claim, suppose on the contrary and up
to symmetry that $z$ has neighbours $w_1\in V_1$ and $w_j\in V_j$ for
some even $j$.  First suppose that $j\in\{2, \ell\}$, say (up to
symmetry) $j=2$.  Pick any $w_h\in V_h$ for $h=3, \ldots, \ell$.  Then
$w_1, \ldots, w_\ell$ induce a hole.  If $w_1, w_2$ are the only
neighbors of $z$ in that hole, then $z$ is a sensitive vertex, a
contradiction.  So $z$ has at least three vertices in that hole and,
by Lemma~\ref{lem:neighbourhood}, $z$ must see exactly one of $w_\ell,
w_3$, say $z$ sees $w_\ell$, and then miss all of $w_3, \ldots,
w_{\ell-1}$ (if $\ell\ge 8$) or miss $w_3, w_5$ (if $\ell=6$).
Repeating this argument for every choice of $w_h$ with $h\neq 1$, we
obtain that $z$ sees all of $V_\ell\cup V_2$ and misses all of
$V_{\ell-1}\cup V_3$.  Since $z$ has a neighbour in $V_1$, $z$ must be
in $D_1$, a contradiction.  Now suppose that $4\le j\le l-2$.  Pick
any $w_h\in V_h$ for $h=2, \ldots, \ell$, $h\neq j$.  Then $w_1,
\ldots, w_\ell$ induce a hole.  By Lemma~\ref{lem:neighbourhood} and
up to symmetry, we must have $\ell=6$, $j=4$, and $z$ must see both
$w_6, w_2$ and miss both $w_{5}, w_3$.  Then repeating this argument
for every choice of $w_h$ with $h\in\{2, 3, 5, 6\}$ implies that $z$
sees all of $V_6\cup V_2$ and misses all of $V_{5} \cup V_3$.  Since
$z$ has a neighbour in $V_4$, $z$ must be in $D_1$, a contradiction.
Thus we have prove the first part of the claim.  To prove the second
part, suppose on the contrary that $z$ sees all of $V_i\cup V_j$ for
some $i\neq j$.  So $i,j$ have the same parity.  If $j=i+2$, then $z$
should be in $D_{i+1} \cup A^*_1 \cup A^*_2$, a contradiction.
Likewise for $j = i-2$.  So $i+4\le j\le i-4$, where indices are taken
modulo $\ell$, which contradicts Lemma~\ref{lem:neighbourhood}.
%Then $\ell\ge 8$, and $z$ has non-neighbours $w_h\in
%V_h$ for each $h\in \{i+1, \ldots, j-1\}$.  But then $z, w_i, w_{i+1},
%\ldots, w_{j-1}, w_j$ induce an even hole of length at least $6$ and
%at most $\ell-2$, a contradiction to the definition of $\ell$.  
Therefore (\ref{zv12}) holds.
\begin{claim}\label{zc12}
%Every vertex of $Z$ misses either all of $C^*_1$ or all of $C^*_2$. 
Each vertex of $Z$ misses all of $C^*_1$ or all of $C^*_2$.
\end{claim}
%
% \subsection*{Proof of claim \ref{zc12} in Theorem \ref{nocycb}}
For suppose that $z$ has neighbours $x\in C^*_1$ and $y\in C^*_2$.  Up
to symmetry there are two cases: (a) $x\in D_1\cup A^*_1$ and $y\in
D_2\cup A^*_2$; and (b) $x\in D_1$ and $y\in D_j$ with $4\le j\le
l-2$.  In either case, by (\ref{zv12}) we can pick vertices $w_h\in
V_h$ for $h = 1, \ldots, \ell$ such that $z$ sees at most one of them.
\\
Consider case (a).  Suppose that $x,y$ are adjacent.  If $z$ sees
$w_4$, then $w_\ell$-$xzy$-$w_3$ is a bull; if $x$ misses $w_4$, then
$w_\ell$-$xyz$-$w_4$ is a second bull, while if $x$ sees $w_4$, then
$w_\ell$-$xzw_4$-$w_3$ is a second bull, a contradiction.  So $z$
misses $w_4$.  Likewise $z$ misses $w_{\ell-1}$.  Suppose that $z$
sees $w_3$.  Then $w_1$-$yzw_3$-$w_4$ is a bull.  Then $x$ misses
$w_1$, for otherwise $w_{\ell-1}$-$w_\ell w_1x$-$z$ is a second bull.
Then $x$ misses $w_4$, for otherwise $w_1$-$yzx$-$w_4$ is a second
bull.  Then $y$ misses $w_2$, for otherwise $w_\ell$-$w_1w_2y$-$z$ is
a second bull.  Then $y$ misses $w_5$, for otherwise $w_5$-$yzx$-$w_2$
is a second bull.  Then $y$ misses $w_4$, for otherwise
$w_2$-$xzy$-$w_4$ is a second bull.  But now vertices $w_1, y, w_3,
w_4, \ldots, w_\ell$ induce a hole, and the neighbors of $z$ in that
hole are $y$ and $w_3$, so $z$ is a sensitive vertex, a contradiction.
So $z$ misses $w_3$.  Similarly $z$ misses $w_\ell$.  Thus
$w_\ell$-$xzy$-$w_3$ is a bull.  If $x$ has a non-neighbour $u_4\in
V_4$ and $y$ has a non-neighbour $u_{\ell-1}\in V_{\ell-1}$ then $x,
y, w_3, u_4, \ldots, u_{\ell-1}, w_\ell$ induce a hole, and the
neighbors of $z$ in that hole are $x$ and $y$, so $z$ is a sensitive
vertex, a contradiction.  So we may assume up to symmetry that $x$
sees every vertex of $V_4$, which, by (\ref{xivi}), implies that
$\ell=6$ and $x$ has a neighbour $u_1\in V_1$.  But then if $z$ sees
$u_1$, then $w_6$-$u_1zy$-$w_3$ is a second bull containing $z$, while
if $z$ misses $u_1$, then $w_5$-$w_6u_1x$-$z$ is a second bull
containing $z$, a contradiction.  Thus $x,y$ are not adjacent.  This
implies that $x\in X_1$ and $y\in X_2$, so $\ell=6$ and $x$ has a
neighbour $u_4\in V_4$ and $y$ has a neighbour $u_5\in V_5$.  Then $z$
sees one of $u_4, u_5$, for otherwise $z, x, u_4, u_5, y$ induce a
$C_5$.  Up to symmetry $z$ sees $u_4$.  Then (\ref{zv12}) implies that
$z$ misses all of $w_1, w_3, w_5$.  If $z$ misses $w_6\in V_6$, then
$w_6$-$xzu_4$-$w_3$ and $w_6$-$xu_4z$-$y$ are two intersecting bulls,
a contradiction.  So $z$ sees $w_6$, and it misses $w_2$.  But then
$w_5$-$w_6zx$-$w_2$ and $w_5$-$u_4zx$-$w_2$ are two intersecting
bulls, a contradiction.  \\
Now consider case (b).  By (\ref{zv12}), we have either $x\in X_1$ or
$y\in X_j$, and so $\ell=6$ and $j=4$.  Then, up to symmetry, $z$ misses
$w_3, w_5, w_6$.  Then $x,y$ are adjacent, for otherwise $z, y, w_5,
w_6, x$ induce a $C_5$.  Then $w_6$-$xzy$-$w_3$ is a bull.  If $z$
misses $w_2$, then $w_2$-$xzy$-$w_5$ is a second bull, while if $z$
sees $w_2$, then $w_6$-$xzw_2$-$w_3$ is a second bull, in either case
a contradiction.  Therefore (\ref{zc12}) holds.

We let $Z^*_1$ (resp.~$Z^*_2$) denote the set of vertices of $Z$ that
have a neighbour in $C^*_2$ (resp.~$C^*_1$).  By Claim~\ref{zc12},
$Z^*_1\cap Z^*_2=\emptyset$, there is no edge between $Z^*_1$ and
$C^*_1$, and there is no edge between $Z^*_2$ and $C^*_2$.
 
\ 

Now, we decompose the whole graph into connected subsets based on a
``hanging'' from $C^*_1$.  Precisely, let us define sets: $$L_1=
C^*_1, \ \ L_2=N(L_1)= C^*_2\cup Z^*_2,\ \ L_{j+1} = N(L_j) - L_{j-1}$$
for any $j\ge 2$, as long as this defines non-empty sets.  The $L_j$'s
will be called the levels of the decomposition.  Note that $Z^*_1
\subseteq L_3$ by (\ref{zc12}).  Level $L_i$ will be called odd or
even according to the parity of $i$.

\

A vertex will be called {\it central\/} if it is in $C^*_1\cup C^*_2$,
and {\it peripheral\/} otherwise.  We will call {\it box\/} any subset
that induces a connected component in any $L_j$.  It is clear that the
whole vertex set of the graph is partitioned into boxes.  By Fact
(\ref{zc12}), a box may contain either only central vertices or only
peripheral vertices.  The boxes will be called central or peripheral
accordingly.  More precisely, by (\ref{nodidj}) and (\ref{noaidi}),
every central box is a subset of some $D_i$ or of some $A^*_i$.  Level
$L_1$ consists of central boxes only.  Level $L_2$ consists of the
central boxes in $C^*_2$, plus the peripheral boxes in $Z^*_2$ (if
any).  The deeper levels consist of peripheral boxes exclusively.
Clearly, Properties (i), (ii), (iii) hold.  We now prove that each box
in $L_j$ satisfies the desired Properties~(iv), (v), (vi) or~(vii) by
induction on $j$.
\begin{claim}\label{bcent}
Every central box satisfies Properties (iv) and (vii).
\end{claim}
%
% \subsection*{Proof of claim \ref{bcent} in Theorem \ref{nocycb}}
For let $B$ be any central box.  We may assume up to symmetry that
$B\subseteq D_3$ or $B\subseteq A^*_1$.  In either case every vertex
of $B$ sees all of $V_2\cup V_4$ and misses all of $V_1\cup V_5$.  We
claim that every $z\in P(B)$ misses all of $V_2\cup V_4$.  For suppose
on the contrary, and up to symmetry, that $z$ sees some $w_2\in V_2$.
There are adjacent vertices $u,v\in B$ such that $z$ sees $u$ and
misses $v$.  Suppose that $z$ sees any $w_1\in V_1$.  Then
(\ref{zv12}) implies $z\in D_1\cup D_2$.  If $z\in D_1$, the edge $zu$
contradicts (\ref{nodidj}) or (\ref{noaidi}).  So $z\in D_2$, and so
$z$ misses all of $V_\ell\cup V_4$.  Then $v_\ell$-$w_1zw_2$-$v$ is a
bull.  If $z$ misses any $w_5\in V_5$, then $z$-$uvv_4$-$w_5$ is a
second bull containing $z$, while if $z$ sees any $w_5\in V_5$,
$w_5$-$zw_1w_2$-$v$ is a second bull containing $z$, a contradiction.
Thus $z$ misses all of $V_1$.  Suppose that $z$ also sees some $w_4\in
V_4$.  Then by symmetry $z$ misses all of $V_5$.  Vertex $z$ cannot
see all of $V_2\cup V_4$, for otherwise $z$ would be in $D_3\cup
A^*_1$, contradicting the fact that $z\in P(B)$.  So, up to symmetry,
we may assume that $z$ has a non-neighbour $w'_2\in V_2$.  But then
$v_1$-$w'_2vu$-$z$ and $w'_2$-$uzw_4$-$v_5$ are two intersecting
bulls, a contradiction.  Thus $z$ misses all of $V_4$.  Then
$v_1$-$w_2zu$-$v_4$ is a bull.  If $z$ misses any $w_5\in V_5$, then
$z$-$uvv_4$-$w_5$ is a second bull containing $z$, while if $z$ sees
$w_5$, then $z, w_5, v_4, v, w_2$ induce a $C_5$, a contradiction.  So
we have proved that every vertex in $V_2\cup V_4$ misses all of
$P(B)$.  Thus it suffices to take auxiliary vertices $a_B = v_2$, $b_B
= v_1$, $c_B=v_4$, $d_B=v_5$ for $B$.  Thus Property (iv) is
established.  To prove (vii), suppose on the contrary that there are
vertices $u,v,w, x, y$ as in the statement of (vii).  Then $v_1$ sees
one of $x,y$, for otherwise $v_1$-$v_2vu$-$x$ and $v_1$-$v_2wv$-$y$
are two intersecting bulls.  Then $v_1$ sees $y$, for otherwise it
sees $x$ and then $v_1, x, y, v, v_2$ induce a $C_5$.  Likewise $v_5$
sees $y$.  But then $v_1$-$yuv$-$w$ and $v_5$-$yuv$-$w$ are two
intersecting bulls, a contradiction.  So (vii) is established.
Therefore (\ref{bcent}) holds.

% About the partials on a central box.  USELESS?
% \\
% Now we prove the last sentence of (\ref{bcent}).  Again let $i=3$.
% Suppose on the contrary that some $w_1\in V_1$ misses a vertex $z\in
% P(B)$ and some $w_5\in V_5$ misses a vertex $z'\in P(B)$.  There are
% adjacent vertices $u,v\in B$ such that $z$ sees $u$ and misses $v$,
% and there are adjacent vertices $u', v'$ such that $z'$ sees $u'$ and
% misses $v'$.  By the first part of this claim, $z$ and $z'$ miss every
% vertex of $V_2\cup V_4$.  Then $w_1$-$v_2vu$-$z$ and
% $w_5$-$v_4v'u'$-$z'$ are two bulls, so they do not intersect, so
% $z\neq z'$ and $u,v,u',v'$ are four distinct vertices.  Then $z$ sees
% $w_5$, for otherwise $w_5$-$v_4vu$-$z$ is a second bull containing
% $z$.  Likewise $z'$ sees $w_1$.  Suppose that $z$ sees one of $u',
% v'$.  Then it sees both, for otherwise $z, u', v', v_2, w_1$ induce a
% second bull containing $z$.  Then $z$ sees $z'$, for otherwise
% $w_5$-$zv'u'$-$z'$ is a second bull.  But then $w_1, v_2, v', z, z'$
% induce a $C_5$, a contradiction.  So $z$ misses both $u', v'$.
% Likewise $z'$ misses both $u,v$.  Then $u$ misses $u'$ and $v'$, for
% otherwise $z$-$uu'v_2$-$w_1$ or $z$-$uv'v_2$-$w_1$ is a second bull
% containing $z$.  Likewise $u'$ misses $u$ and $v$.  But then $z$-$u v
% v_2$-$u'$ is a second bull containing $z$, a contradiction. 

The preceding claims imply that all boxes in $L_1$ and all boxes in
$L_2-Z^*_2$ satisfy Properties~(iv) and (vii).  Now we consider the
peripheral boxes, which are the boxes in $Z^*_2$ and in
$L_j$ for any $j\ge 3$.  First we consider the boxes in $Z^*_2$.
\begin{claim}\label{rab}
%Given non-adjacent vertices $a,b$, both in $C^*_1\cup Z^*_1$ or in
%$C^*_2\cup Z^*_2$, 
Given non-adjacent vertices $a,b$, both in $C^*_1$, 
or both in $C^*_2\cup Z^*_2$, or both in $Z^*_1$,
there exists a chordless even path $R_{ab}$ whose
interior vertices are alternately in $L_1$ and $L_2$.
%$V^*_1$ and $V^*_2$.
\end{claim} 
%
% \subsection*{Proof of claim \ref{rab} in Theorem \ref{nocycb}}
 For suppose first that $a,b$ are both in $L_{1} = C^*_1$.
%For suppose up to symmetry that $a,b$ are both in $C^*_1\cup Z^*_1$.  
By the definition of the $D_i$'s and $A^*_j$'s,
%and $Z^*_j$'s, 
every vertex in any such set is adjacent to some vertex of $V^*$.
Thus there is a path from $a$ to $b$ whose interior vertices are
alternately in even $V_i$'s and odd $V_i$'s and no two interior
vertices are in the same $V_i$.  Take a shortest such path $R=
a$-$v_h$-$\cdots$-$v_k$-$b$.  Clearly, $R$ has even length.  Then
(\ref{nodidj}) and (\ref{noaidi})
%(\ref{zc12}) 
imply that any chord of $R$ must be of the type $av_i$ for some even
$i>h$ or $bv_j$ for some even $j<k$, and so
$a$-$v_i$-$\cdots$-$v_k$-$b$ or $a$-$v_h$-$\cdots$-$v_j$-$b$ is a path
with the same properties and shorter than $R$, a contradiction.  Now
suppose that $a,b$ are both in $L_2 = C^*_2\cup Z^*_2$.  Let $a'$ be a
neighbour of $a$ in $L_{1}$ and $b'$ be a neighbour of $b$ in $L_{1}$.
If $a,b$ have a common such neighbour, then we can take $a'=b'$ and
$R_{ab} = a$-$a'$-$b$.  In the remaining case, we may assume that $a$
misses $b'$ and $b$ misses $a'$.  If $a', b'$ are adjacent, then they
lie in one box in $L_{1}$, for which Property (iv) is already proved,
and then $a, a', b, b'$ violate Lemma~\ref{lem:moreitems}.  So $a',b'$
are not adjacent vertices in $L_{1} =C^*_1$ and there exists a path
$R_{a'b'}$ with the desired properties.  Then the path
$a$-$R_{a'b'}$-$b$ has even length and its interior vertices are in
$L_1\cup L_2$ and alternately in odd and even $L_i$'s.  If this path
has any chord, then it must be incident with $a$ or $b$, and then
(\ref{zc12}) implies that we can find a shorter subpath with the same
properties.  Finally, suppose that $a,b$ are both in $Z^*_1$.  Note
that the definition of the levels implies that the set $Z^*_1$ is
contained in $L_{3}$.  So, by considering a neighbour $a'$ of $a$ in
$L_{2}$ and a neighbour $b'$ of $b$ in $L_{2}$, and by applying an
analogous argument we obtain a path with the desired properties.
Therefore (\ref{rab}) holds.
\begin{claim}\label{bz2}
Every box in $Z^*_2$ and every box in $Z^*_1$ satisfy 
Properties (v), (vi), (vii).
\end{claim}
%
% \subsection*{Proof of claim \ref{bz2} in Theorem \ref{nocycb}}
For let $B$ be any box in $Z^*_2$.  First let us prove the assertion
that, for every subset $C\subseteq B$ that induces a connected
subgraph, there is a vertex of $L_1$ that sees all vertices of $C$.
We prove the assertion by induction on $|C|$.  If $|C|=1$ the
assertion holds by the definition of $B$.  Now suppose that the
assertion holds for any $C$ of size at most $k$, and let $C$ have size
$k+1\ge 2$.  Let $c_1$-$\cdots$-$c_h$ be a longest chordless path in
$C$.  Thus $C-c_1$ and $C-c_h$ are connected and, by the induction
hypothesis, there is a vertex $u\in L_1$ that sees all of $C-c_1$, and
there is a vertex $v\in L_1$ that sees all of $C-c_h$.  If $u$ sees
$c_1$, or $v$ sees $c_h$, then we are done.  So let us assume that $u$
misses $c_1$ and $v$ misses $c_h$.  Note that for each $a$ in $L_1$
there is a vertex $a'$ that sees $a$ and misses all of $B$; indeed,
$a$ is in $D_i\cup A^*_1$ for some odd $i$, and so any vertex in
$V_{i+1}$ can play the role of $a'$.  In particular we can consider
vertices $u'$ and $v'$.  If $h\ge 6$, then $c_1$-$c_2c_3u$-$c_5$ and
$c_1$-$c_2c_3u$-$c_6$ are two intersecting bulls, a contradiction.  If
$3\le h\le 5$, then $c_1$-$c_2c_3u$-$u'$ and
$c_h$-$c_{h-1}c_{h-2}v$-$v'$ are two intersecting bulls, a
contradiction.  So $h=2$.  This means that $C$ is a clique.  Suppose
that $u$ misses $v$.  Consider any path $R_{uv}=r_1$-$\cdots$-$r_p$
given by~(\ref{rab}), with $p$ odd, $r_1=u$, $r_p=v$.  Then
$v$-$c_1$-$c_2$-$u$-$R_{uv}$-$v$ is an odd cycle of length at least
five, so it must contains a triangle, for otherwise it contains an odd
hole.  Note that $c_1$ and $c_2$ do not see two consecutive vertices
on the path $R_{uv}$, since they are in $Z^*_2$ and by (\ref{zc12}).
So, in order to have a triangle, there must be a vertex $r_j$ that
sees both $c_1, c_2$, and so $r_j\in L_1$, and so $3\le j\le p-2$.
But then $u$-$c_2c_1r_j$-$r_{j+1}$ and $v$-$c_1c_2r_j$-$r_{j-1}$ are
two intersecting bulls, a contradiction.  Thus $u,v$ are adjacent, and
so they lie in one box $U$ of $L_1$, and $c_1, c_2\in P(U)$.  Up to
symmetry, we may assume that $U\subseteq D_3\cup A^*_1$ and so, as
proved in (\ref{bcent}), $v_2$ sees all of $U$ and misses all of
$P(U)$ and $v_1$ misses all of $U$.  Then $v_1$ sees one of $c_1,
c_2$, for otherwise $v_1$-$v_2uv$-$c_1$ and $v_1$-$v_2vu$-$c_2$ are
two intersecting bulls.  Then $v_1$ sees both $c_1, c_2$, for if it
sees only one, say $c_1$, then $v_1, v_2, u, c_2, c_1$ induce a $C_5$.
Then $v_1$ sees every $z\in C-\{c_1, c_2\}$, for otherwise $v_1, v_2,
u, z, c_1$ induce a $C_5$ (recall that $C$ is a clique and $u$ sees
all of $C-c_1$).  Thus we have proved the assertion.  Applying it to
$C=B$, we obtain that some vertex $a$ of $L_1$ sees all of $B$.  \\
Up to symmetry we may assume that $a\in D_3\cup A^*_1$.  Now we claim
that $a$ misses every vertex of $P(B)$.  For suppose on the contrary
that $a$ sees some $x\in P(B)$.  There are adjacent vertices $u,v\in
B$ such that $x$ sees $u$ and misses $v$.  Since $x$ sees $a$, we have
$x\in D_2\cup D_4\cup D_6\cup A^*_2\cup Z^*_2$.  In fact we do not
have $x\in Z^*_2$, for otherwise $x$ should be in $B$.  Also we do not
have $x\in D_6$, for otherwise $u$ would contradict (\ref{zc12}). %
Thus, up to symmetry, we have $x\in D_2\cup A^*_2$.  Since $u,v\in
Z^*_2$, they both miss all of $V_2\cup V_4$.  Up to symmetry we may
assume that $u$ misses some $w_1\in V_1$.  Then $w_1$-$xua$-$v_4$ is a
bull.  Then $v$ misses $w_1$, for otherwise $w_1$-$vua$-$v_4$ is a
second bull containing $a$.  Suppose that $x$ sees every $v_5\in V_5$.
If $x$ sees any $w_2\in V_2$, then $v$ sees $v_5$, for otherwise
$v_5$-$xv_2a$-$v$ is a second bull containing $a$, and then $u$ sees
every $v_5$, for otherwise $v_5$-$vua$-$v_2$ is a second bull
containing $a$.  On the other hand if $x$ misses any $w_2\in V_2$,
then $u$ sees $v_5$, for otherwise $v_2$-$aux$-$v_5$ is a second bull.
In either case $u$ sees all of $V_5$.  Since $u\in Z^*_2$, $u$ must
miss some $w_3\in V_3$.  But then $v_6$-$v_5ux$-$w_3$ is a second bull
containing $x$, a contradiction.  %
Thus $x$ has a non-neighbour $w_5\in V_5$.  Suppose $a$ has a
non-neighbour $w_6\in V_6$.  If $\ell\ge 8$, let us pick any $w_i\in
V_i$ for $i=7, \ldots, \ell$; then $w_1, v_2, a, v_4, w_5, \ldots,
w_\ell$ induce a hole in $G-u$, and $w_1$-$x$-$a$-$v_4$-$w_5$-$w_6$ is
a $P_6$ or $C_6$ in $G-u$ such that $u$ is adjacent to $x,a$ and not
to $w_1, v_4, w_5, w_6$, so $u$ is a sensitive vertex, a
contradiction.  So $a$ sees all of $V_6$, and so $a\notin V_3$.
Consider any $v_3\in V_3$.  The same argument as for $a$ implies that
$v_3$ misses one of $u,v$.  Then $v_3$ misses $a$, for otherwise
$w_1$-$v_2v_3a$-$u$ or $w_1$-$v_2v_3a$-$v$ is a second bull.  Then
$v_3$ sees $u$, for otherwise $v_3$-$xua$-$v_6$ is a second bull.  So
$v_3$ misses $v$.  But then $v_3$-$uva$-$v_6$ is a second bull
containing $a$, a contradiction.  Therefore $a$ misses all of $P(B)$.
So $a$ can play the role of $a_B$, and any vertex in $V_2$ can play
the role of $b_B$.  Thus Property~(v) is established.  \\
% Now we claim that some vertex in $V_1\cup V_5$ misses all of $B$.  To
% establish this, let us prove the assertion that, for every subset
% $C\subseteq B$ that induces a connected subgraph, there is a vertex of
% $V_1\cup V_5$ that misses all of $C$.  We prove this by induction on
% $|C|$.  If $|C|=1$, the fact follows from (\ref{zc12}).  So suppose
% that the assertion holds for every set of size $k$ and let $|C|=k+1$.
% Let $c_1$-$\cdots$-$c_h$ be a longest chordless path in $C$, with
% $h\ge 2$.  So $C-c_1$ and $C-c_h$ are connected, and by induction
% there are vertices $x,y\in V_1\cup V_5$ such that $x$ misses all of
% $C-c_1$ and $y$ misses all of $C-c_h$.  We may assume that $x$ sees
% $c_1$ and $y$ sees $c_h$, for otherwise we are done.  Then one of
% $x$-$c_1c_2a$-$v_2$ or $x$-$c_1c_2a$-$v_4$ is a bull, and similarly
% one of $y$-$c_hc_{h-1}a$-$v_2$ or $x$-$c_hc_{h-1}a$-$v_4$ is a bull,
% so there are two bulls containing $a$, a contradiction.  So the
% assertion holds.  Applying it to $C=B$ yields, up to symmetry, a
% vertex $w_1\in V_1$ that misses all of $B$.  \\ 
In order to prove Property~(vi), suppose on the contrary that there
are two adjacent vertices $u,v\in B$ and two vertices $x, y\in P(B)$
such that $x$ sees $u$ and misses $v$ and $y$ sees $v$ and misses $u$.
Since $B$ satisfies~(v), by Lemma~\ref{lem:moreitems}, $x$ sees $y$.
If $x$ misses any $w_2\in V_2$, then $y$ also misses $w_2$, for
otherwise $w_2, a, u, x, y$ induce a $C_5$; but then $w_2$-$avu$-$x$
and $w_2$-$auv$-$y$ are two intersecting bulls, a contradiction.  Thus
$x$ sees all of $V_2$.  Similarly $x$ sees all of $V_4$.  But then $x$
should be in $D_3\cup A^*_1$, a contradiction.  Thus Property~(vi) is
established.  \\
In order to prove (vii), suppose on the contrary that there are
vertices $u, v, w, x, y$ as in the statement of (vii).  Then $v_2$
sees one of $x,y$, for otherwise $v_2$-$avu$-$x$ and $v_2$-$awv$-$y$
are two intersecting bulls.  Then $v_2$ sees $y$, for otherwise it
sees only $x$, and then $v_2, x, y, v, a$ induce a $C_5$.  Likewise
$v_4$ sees $y$.  But then $v_2$-$yuv$-$w$ and $v_4$-$yuv$-$w$ are two
intersecting bulls.  So (vii) is established.

An analogous argument establishes the properties for a box in $Z^*_1$.
Therefore (\ref{bz2}) holds.

Now we consider the boxes in $L_j$ for $j\ge 3$ that are not in $Z^*_1$.  
\begin{claim}\label{rabdeep}
For any $j\ge 1$, given non adjacent vertices $a, b$ in $L_j$, there
exists a chordless even path $R_{ab}$ whose interior vertices are in
$L_1\cup\cdots\cup L_{\max\{2, j-1\}}$ and alternately in odd and even
$L_i$'s.  Moreover, every box in $L_j$ with $j\ge 3$ satisfies
Properties~(v), (vi) and~(vii).
\end{claim}
%
% \subsection*{Proof of claim \ref{rabdeep} in Theorem \ref{nocycb}}
We prove this claim by induction on $j$.  For $j=1$, or $j=2$, or
$j=3$ and $a,b$ both in $Z^*_1$, this is Claims~(\ref{bcent}),
(\ref{rab}), (\ref{bz2}).  Now suppose $j\ge 3$ and $a,b$ not both in
$Z^*_1$.  Let $a'$ be a neighbour of $a$ in $L_{j-1}$ and $b'$ be a
neighbour of $b$ in $L_{j-1}$.  If $a,b$ have a common such neighbour,
then we can take $a'=b'$ and $R_{ab} = a$-$a'$-$b$.  In the remaining
case, we may assume that $a$ misses $b'$ and $b$ misses $a'$.  If $a',
b'$ are adjacent, then they lie in one box in $L_{j-1}$, for which
Property (iv) or (v) is already proved, and then $a, a', b, b'$
violate Lemma~\ref{lem:moreitems}.  So $a',b'$ are not adjacent.  By
the induction hypothesis, there exists a path $R_{a'b'}$ with the
desired properties.  Then the path $a$-$R_{a'b'}$-$b$ has even length
and its interior vertices are in $L_1\cup\cdots\cup L_{j-1}$ and
alternately in odd and even $L_i$'s.  If this path has any chord, then
it must be incident with $a$ or $b$, and then (\ref{zc12}) implies
that we can find a shorter subpath with the same properties.
%Therefore (\ref{rabdeep}) holds.
%
%
%\begin{claim}\label{bz3}
%Every box in $L_j$ with $j\ge 3$ satisfies Properties~(v), (vi)
%and~(vii).
%\end{claim}
%
% \subsection*{Proof of claim \ref{bz3} in Theorem \ref{nocycb}}

Now the proof of Properties~(v), (vi) and~(vii)
is rather similar to the proof of (\ref{bz2}).
First we prove that some vertex of $L_{j-1}$ sees all of $B$, with the
following changes: instead of (\ref{rab}), use the chordless even path
$R_{ab}$ given by induction;
%(\ref{rabdeep}); 
for every vertex $a$ in $L_{j-1}$, there is a neighbour $a'$ of $a$ in
$L_{j-2}$ that misses all of $B$; when (\ref{bcent}) is invoked to
deal with the $C_4$ induced by $u,v,c_1, c_2$ with $u,v\in U$, we can
still invoke (\ref{bcent}) if $U$ is a central box, and we can invoke
Property~(vi) when $U$ is a peripheral box, since that property holds
for $U$ by the induction hypothesis on $j$.  Thus there is a vertex
$a\in L_{j-1}$ that sees all of $B$.  If $a$ is in a central box, the
rest of the proof is completely the same, with subscripts shifted by
$1$.  There remains to deal with the case when $a$ is in $Z$.  \\
We prove that $a$ misses all of $P(B)$.  %
% \footnote{%%%%%%
% This is where we need the absence of pericyclic bull.  For example,
% consider vertices $z, x, a, u, v$ in $Z$ such that $z$ sees only
% $v_1$, vertices $x,a,u,v$ induce a diamond (with non-edge $xv$) and
% miss all of $v_1, \ldots, v_\ell$, and $z$ sees $x,a$.  This graph has
% only one bull.  Then $\{u,v\}$ is a peripheral box, but it does not
% have a vertex $a_B$ that sees all of $B$ and misses all of $P(B)$.  In
% this example, the bull is pericyclic.  We use the box partition in the
% proof of Theorem~\ref{mbpto}, so we need these properties!!  The
% hypothesis `No pericyclic bull' forces $a$ to miss every vertex of
% $P(B)$, so (v) holds, as proved here.  How about (vi) and (vii)?  I
% did not have time to finish\ldots}
%%%.  
Suppose that $a$ sees a vertex $x\in P(B)$.  So there are adjacent
vertices $u,v\in B$ such that $x$ sees $u$ and misses $v$.  Since $x$
is not in $B$, it is not in $L_j$; and since it sees $a$ and $u$, it
must be in $L_{j-1}$.  So $a,x$ are in a box $U\subseteq L_{j-1}$, and
this box has auxiliary vertices $a_U, b_U$ by the induction hypothesis
on $j$, with $a_U\in L_{j-2}$.  Suppose that $a_U$ is a central
vertex, say $a_U\subseteq D_3\cup A^*_1$ with $a\in Z^*_2$ (the case
$a_U\in D_2\cup A^*_2$ is similar).  So $U\subseteq Z^*_2$, and so
$a,x$ miss all of $V^*_2$.  Vertices $u,v$ miss all of $V^*_1$ since
they are in $L_3$.  Since $u\in Z$, by (\ref{zv12}) it has a
non-neighbour in $V_2\cup V_4$, say $u$ misses $w_4\in V_4$.  Then $v$
misses $w_4$, for otherwise $w_4, v, u, x, a_U$ induce a $C_5$.  So
$v$-$axa_U$-$w_4$ is one bull.  Then $v$ sees every $v_2\in V_2$, for
otherwise $v$-$axa_U$-$v_2$ is a second bull; and $u$ sees every
$v_2$, for otherwise $v_2, v, u, x, a_U$ induce a $C_5$.  By
(\ref{zv12}), $x$ misses some $w\in V_1\cup V_3$.  But then
$w$-$v_2vu$-$x$ is a second bull.  So $a_U$ is not a central vertex,
and so $j\ge 4$.  By the definition of the levels, there is a shortest
path $p_1$-$\cdots$-$p_r$ such that $p_r=a_U$ and $p_1$ is in
$Z^*_1\cup Z^*_2$ (and so every vertex of $P\setminus p_1$ has no
neighbour in $C^*_1\cup C^*_2$).  By (\ref{zv12}), there are vertices
$w_i\in V_i$, $i=1,\ldots, \ell$, such that $p_1$ sees exactly one of
them.  If $j\ge 5$, then the subgraph of $G-x$ induced by vertices
$w_1, \ldots, w_\ell, p_1, \ldots, p_R=a_U, a, v$ contains a hole and
a $P_6$ $v$-$a$-$a_U$-$\cdots$ such that $x$ sees $a_U, a$ and misses
the other four vertices of the $P_6$, so $x$ is a sensitive vertex, a
contradiction.  So $j=4$, and so $a_U=p_1$.  Since $a,x$ are in $L_3$
they miss all of $w_1, w_3, \ldots, w_{\ell-1}$.  Suppose that some
$w_j$ with even $j$ sees $x$; then it sees $a$, for otherwise
$w_j$-$xa_Ua$-$v$ is a second bull; but then $w_{j-1}$-$w_jxa$-$v$ is
a second bull, a contradiction.  So $x$ misses every $w_i$.  But then
the subgraph of $G-x$ induced by vertices $w_1, \ldots, w_\ell, a_U,
a, v$ contains a hole and a $P_6$ $v$-$a$-$a_U$-$w_1$-$w_2$-$w_3$ such
that $x$ sees $a_U, a$ and misses the other four vertices of the
$P_6$, so $x$ is a sensitive vertex, a contradiction.  Thus we have
proved that $a$ misses every vertex of $P(B)$.  Since $j\ge 3$, $a$
has a neighbour $a'$ in $L_{j-2}$, and so $a$ and $a'$ can play the
role of $a_B$ and $b_B$, and Property (v) is established.  \\
There remains to prove (vi) and (vii).  In order to prove
Property~(vi), suppose on the contrary that there are two adjacent
vertices $u,v\in B$ and two vertices $x, y\in P(B)$ such that $x$ sees
$u$ and misses $v$ and $y$ sees $v$ and misses $u$.  Since $B$
satisfies~(v), by Lemma~\ref{lem:moreitems}, $x$ sees $y$.  If $a'$
misses both $x,y$, then $a'$-$auv$-$y$ and $a'$-$avu$-$x$ are two
intersecting bulls; and if $a'$ sees only one of $x,y$, then $a', a,
x, y$ and one of $u,v$ induce a $C_5$.  So $a'$ sees both $x,y$.  Thus
$x,y$ are in one box $B'$, which is in $L_{j-1}$.  If $j\ge 4$, then
$a'$ has a neighbour $a''$ in $L_{j-3}$, and $a''$-$a'xy$-$v$ and
$a''$-$a'yx$-$u$ are two intersecting bulls.  So $j=3$, and $x,y$ are
in $L_2$.  If the box $B'$ that contains $x,y$ is peripheral, then the
situation contradicts the fact that $B'$ satisfies (vi), which was
proved in (\ref{bz2}).  So $B'$ is a central box.  Up to symmetry we
may assume that $B'\subset V_4\cup X_4\cup A^*_2$, and so, as in the
proof of (\ref{bcent}), we know that $v_3$ and $v_5$ see all of $B'$
and every vertex in $V_2\cup V_6$ misses all of $B'$.  In consequence
every $w$ in $V_2$ sees both $u,v$ (for if it misses both then
$w$-$v_3xy$-$u$ and $w$-$v_3yx$-$u$ are two intersecting bulls, and if
it sees only one then $w, v_3, u, v$ and one of $x,y$ induce a $C_5$);
and similarly every $w$ in $V_6$ sees both $u,v$.  But then the fact
that $u,v$ see all of $V_2\cup V_6$ contradicts the definition of the
$X_i$'s, $A^*_j$'s and $Z$.  Thus Property~(vi) is established.  \\
In order to prove (vii), suppose on the contrary that there are
vertices $u, v, w, x, y$ as in the statement of (vii).  If $a'$ misses
both $x,y$, then $a'$-$avu$-$x$ and $a'$-$awv$-$y$ are two bulls.  If
$a'$ sees $x$ and not $y$, then $a', x, y, v, a$ induce a $C_5$.  So
$a'$ sees $y$, and $a'$-$yuv$-$w$ is one bull.  Then $a'$ sees $x$,
for otherwise $a'$-$avu$-$x$ is a second bull.  So $x,y$ are in one
box $B'$ in $L_{j-1}$.  If $j\ge 4$, then $a'$ has a neighbour $a''$
in $L_{j-3}$, and $a''$-$a'xy$-$v$ is a second bull, a contradiction.
So $j=3$.  Since $a'$-$yuv$-$w$ is a bull for each neighbour $a'$ of
$a$ in $L_1$, this $a'$ must be unique, so $a$ is a peripheral vertex.
Since $a'$ is in $L_1$, we may assume up to symmetry that it is in
$V_3\cup X_3\cup A^*_1$.  Then we may assume that $x,y$ are different
from and not adjacent to $v_2$ (else replace $v_2$ by $v_4$).  Then
$v_2$ sees $v$, for otherwise $v_2$-$a'xy$-$v$ is a second bull; $v_2$
sees $u$, for otherwise $v_2, a', x, u, v$ induce a $C_5$; and $v_2$
sees $w$, for otherwise $a'$-$v_2uv$-$w$ is a second bull.  Then $v_1$
sees $x$, for otherwise $v_1$-$v_2vu$-$x$ is a second bull; and $v_1$
sees $y$, for otherwise $v_1$-$v_2wv$-$y$ is a second bull.  But then
$v_1$-$yuv$-$w$ is a second bull.  So (vii) is established.
%Therefore (\ref{bz3}) holds.
Therefore (\ref{rabdeep}) holds.

This completes the proof of Theorem~\ref{nocycb}.  $\Box$

\begin{lemma}\label{lem:sensto}
Let $G$ be a bull-reducible graph that contains no $C_5$, no wheel and
no spiked bull.  Suppose that $G$ has a sensitive vertex $x$, and that
$G-x$ is transitively orientable.  Then $G$ is transitively
orientable.
\end{lemma}
%
% \subsection*{Proof of Lemma \ref{lem:sensto}}
\emph{Proof.} Since $x$ is a sensitive vertex, there exist vertices
$u_1, \ldots, u_6$ of $G-x$ with edges $u_iu_{i+1}$ ($i=1, \ldots, 5$)
and possibly $u_1u_6$, such that $x$ is adjacent to $u_2$ and $u_3$
and not to $u_1, u_4, u_5, u_6$.  Note that $u_1$-$u_2xu_3$-$u_4$ is
one bull, henceforth the ``first bull''.  (Every second bull we will
find will obviously intersect the first one.)  Define sets:
\begin{eqnarray*}
A & = & \{v\in V(G)\mid \mbox{ $v$ sees $x, u_2$ and misses $u_3,
u_5$}\}, \\
B & =& \{v\in V(G)\mid \mbox{ $v$ sees $x, u_3$ and misses $u_1, u_2, 
u_4, u_6$}\}.
\end{eqnarray*}
We first claim that $N(x)=\{u_2, u_3\}\cup A\cup B$.  To prove this,
consider any neighbour $v$ of $x$ different from $u_2, u_3$.  Suppose
that $v$ misses both $u_2, u_3$.  Then $v$ sees $u_1$, for otherwise
$u_1$-$u_2u_3x$-$v$ is a second bull; and similarly $v$ sees $u_4$;
but then $v, u_1, \ldots, u_4$ induce a $C_5$.  So $v$ sees at least
one of $u_2, u_3$.  Suppose that $v$ sees both $u_2, u_3$.  If $v$
sees $u_4$, then it sees $u_1$, for otherwise $u_1$-$u_2xv$-$u_4$ is a
second bull; $v$ sees $u_5$, for otherwise $u_1$-$vu_3u_4$-$u_5$ is a
second bull; and $v$ sees $u_6$, for otherwise $u_2$-$vu_4u_5$-$u_6$
is a second bull; but then, if $u_1u_6$ is not an edge, then $u_1,
\ldots, u_4, x, v, u_6$ induce a spiked bull, and if $u_1u_6$ is an
edge, then $v, u_1, \ldots, u_6$ induce a wheel.  So $v$ misses $u_4$.
Then $v$ misses $u_1$, for otherwise $u_1$-$vxu_3$-$u_4$ is a second
bull.  But then $u_1$-$u_2vu_3$-$u_4$ is a second bull.  Thus $v$ sees
exactly one of $u_2, u_3$.  \\
Now suppose that $v$ sees $u_2$ and misses $u_3$.  Then $v$ misses
$u_5$, for otherwise either $u_1$-$u_2xv$-$u_5$ (if $v$ misses $u_1$)
or $u_3$-$u_2u_1v$-$u_5$ (if $v$ sees $u_1$) is a second bull.  Thus
$v$ is in $A$.  \\
Now, suppose that $v$ sees $u_3$ and misses $u_2$.  If $v$ sees $u_4$,
then $v$ sees $u_5$, for othewise $u_2$-$u_3vu_4$-$u_5$ is a second
bull; and $v$ sees $u_6$, for otherwise $x$-$vu_4u_5$-$u_6$ is a
second bull; but then $u_2$-$u_3u_4v$-$u_6$ is a second bull.  So $v$
misses $u_4$.  Then $v$ misses $u_1$, for otherwise
$u_1$-$vxu_3$-$u_4$ is a second bull; and similarly $v$ misses $u_6$.
Thus $v$ is in $B$.  So we have proved the claim that $N(x)=\{u_2,
u_3\}\cup A\cup B$.

Next, we claim that every vertex in $A$ sees every vertex in $B$.  For
suppose on the contrary that there are non-adjacent vertices $a\in A$,
$b\in B$.  Then $a$ sees $u_4$, for otherwise $a$-$xbu_3$-$u_4$ is a
second bull; but then $b$-$xu_2a$-$u_4$ is a second bull, a
contradiction.  In summary, the two sets $A\cup\{u_3\}$ and
$B\cup\{u_2\}$ form a partition of $N(x)$ and are completely adjacent
to each other.

Let $U_2$ be the set of vertices that see $u_1$ and $u_3$ and miss
$u_4, u_5, u_6$.  Note that $x$ has only one neighbour (which is
$u_2$) in $U_2$, because for any such vertex $w$ there is a bull
$u_1$-$wxu_3$-$u_4$.  Let $D$ be the component of $U_2$ that contains
$u_2$.  Let $N_2=U_2\cap N(u_2)$ and $M_2=U_2- (N_2\cup\{u_2\})$.
Then:
\begin{equation}\begin{minipage}{0.85\linewidth}\label{n2}
Every vertex of $N_2$ sees every vertex of $M_2$, and (consequently)
either $D=\{u_2\}$ or $D=U_2$.
\end{minipage}\end{equation}
For consider any $v\in N_2$ and $w\in M_2$.  Then $v$ sees $w$, for
otherwise $x$-$u_2vu_1$-$w$ is a bull. Therefore (\ref{n2}) holds.
 
\begin{equation}\begin{minipage}{0.85\linewidth}\label{pd}
If $P(D)-x\neq\emptyset$, then $M_2=\emptyset$ and every vertex
$z$ of $P(D)-x$ satisfies one of the following: \\
(a) $z$ sees all of $\{x, u_1, u_3, u_5\}\cup N_2$ and none of $\{u_2,
u_4, u_6\}$; \\
(b) $z$ sees all of $\{u_2, u_4\}$ and none of $\{x, u_1, u_3\}\cup
N_2$.
\end{minipage}\end{equation}
To prove this, suppose that $P(D)-x\neq\emptyset$ 
and let $z$ be any vertex
in $P(D)-x$.  So there are vertices $u,v$ in $D$ such that $z$ sees
$u$ and misses $v$.  By (\ref{n2}), we have $D=U_2$.  So $z$ is not in
$U_2$.  First suppose that $z$ sees both $u_1, u_3$.  If $z$ sees
$u_4$, then it sees $u_5$ (for otherwise $u_1$-$zu_3u_4$-$u_5$ is a
second bull), and it sees $u_6$ (for otherwise $u$-$zu_4u_5$-$u_6$ is
a second bull); but then $v$-$u_3u_4z$-$u_6$ is a second bull.  So $z$
misses $u_4$.  Then $z$ misses $u_6$, for otherwise
$u_4$-$u_3uz$-$u_6$ is a second bull.  Then $z$ sees $u_5$, for
otherwise $z$ should be in $U_2$.  If $x$ misses $z$, then $x$ sees
$u$, for otherwise $x$-$u_3uz$-$u_5$ is a second bull; but then
$x$-$uu_1z$-$u_5$ is a second bull.  So $x$ sees $z$.  Then $x$ sees
$v$, for otherwise $v$-$u_3xz$-$u_5$ is a second bull.  Thus $v=u_2$,
and $u\in N_2$.  Then $z$ sees every $u'\in N_2$, for otherwise
$u'$-$u_3xz$-$u_5$ is a second bull.  If there is any $y\in M_2$, then
$y$ sees $u$ by (\ref{n2}), and $z$ misses $y$, for otherwise
$u_2$-$uyz$-$u_5$ is a second bull; but then $y$-$u_3xz$-$u_5$ is a
second bull.  So $M_2=\emptyset$ and $z$ satisfies (a).  \\
Now suppose that $z$ sees $u_3$ and misses $u_1$.  Then $z$ sees
$u_4$, for otherwise $u_1$-$uzu_3$-$u_4$ is a second bull; and $z$
sees $u_5$, for otherwise $v$-$u_3zu_4$-$u_5$ is a second bull; but
then $u_1$-$uu_3z$-$u_5$ is a second bull.  Therefore, $z$ misses
$u_3$.  Thus $z$ sees $u_4$, for otherwise $z$-$uvu_3$-$u_4$ is a
second bull; and $z$ misses $u_1$, for otherwise $z, u_1, v, u_3, u_4$
induce a $C_5$.  If $z$ sees $x$, then we must have $z\in A$, but then
$u_1$-$u_2xz$-$u_4$ is a second bull.  So $z$ misses $x$.  Then $x$
misses $v$, for otherwise $x$-$vu_1u$-$z$ is a second bull; and $x$
sees $u$, for otherwise $x$-$u_3vu$-$z$ is a second bull.  So $u=u_2$.
Then $z$ misses every $v'$ in $N_2$, for otherwise $x$-$u_2v'z$-$u_4$
is a second bull.  If there is any $y\in M_2$, then $y$ sees $v$ by
(\ref{n2}), and $z$ sees $y$, for otherwise $z$-$u_2xu_3$-$y$ is a
second bull; but then $x$-$u_3vy$-$z$ is a second bull.  Thus
$M_2=\emptyset$ and $z$ satisfies (b).  Therefore (\ref{pd}) holds.

By the hypothesis, there is a transitive orientation of $G-x$.  In
that orientation, we write $u\rightarrow v$ whenever the edge $uv$
exists in $G-x$ and is oriented from $u$ to $v$; and for disjoint sets
$Y, Z\subset V(G)$, we also write $Y\rightarrow Z$ if $y\rightarrow z$
holds for all $y\in Y$ and $z\in Z$.  In the transitive orientation,
we may assume up to symmetry that $u_i\rightarrow u_{i+1}$ for $i=1,
3, 5$ and $u_i\rightarrow u_{i-1}$ for $i=3, 5$.  Then the
transitivity implies $A\rightarrow u_2$, $u_3\rightarrow B$, and
$A\rightarrow B$.  We claim that:

\begin{equation}\begin{minipage}{0.85\linewidth}\label{inu2}
We may assume that every edge $u_2v$ with $v\in U_2$ satisfies 
$v\rightarrow u_2$.  
\end{minipage}\end{equation}
To prove this, first suppose that $P(D)-x=\emptyset$.  So $U_2$ is a
homogeneous set in $G-x$.  Moreover, by (\ref{n2}), every vertex of
$\{u_2\}\cup M_2$ sees every vertex of $N_2$.  So we can reorient the
edges between these two sets in such a way that $N_2\rightarrow
\{u_2\}\cup M_2$.  Then it is easy to see that the modified
orientation is transitive.  Now suppose that $P(D)-x\neq\emptyset$.
So, by (\ref{pd}), we have $U_2=\{u_2\}\cup N_2$.  Let $z$ be any
vertex in $P(D)-x$.  Suppose that $z$ satisfies (a) of (\ref{pd}).
Then the transitivity implies $\{u_1, u_3, u_5\}\rightarrow z$, and,
consequently, $v\rightarrow z$ for every $v\in N_2$, and $v\rightarrow
u_2$ as well.  Thus we have the desired property.  Finally, suppose
that $z$ satisfies (b).  Then the transitivity implies $z\rightarrow
\{u_2, u_4\}$ and consequently $v\rightarrow u_2$ for every $v\in
N_2$.  Thus we also have the desired property.  Therefore (\ref{inu2})
holds.

Let us extend this transitive orientation of $G-x$ to an orientation
of $G$ by setting $a\rightarrow x$ for every $a\in A\cup\{u_3\}$ and
$x\rightarrow b$ for every $b\in B\cup\{u_2\}$.  We claim that this is
a transitive orientation of $G$.  Note that there is no circuit in
$G$, for if a set $S$ of vertices induces a circuit, then $S$ must
contain $x$, and then (since $N(x)=\{u_2, u_3\}\cup A\cup B$ and
$A\cup\{u_3\} \rightarrow B\cup\{u_2\}$) the set $S-x$ would induce a
circuit in $G-x$.  Now suppose that there is a triple $r, s, t$ with
$r\rightarrow s\rightarrow t$ and $r,t$ are not adjacent.  Clearly $x$
is one of $r,s,t$, since the orientation is transitive in $G-x$.  If
$x=s$, then $r$ is in $A\cup\{u_3\}$ and $t$ is in $B\cup\{u_2\}$, but
then we have $r\rightarrow t$ as mentioned above.  So $x\neq s$.  This
leads to the following four cases.

{\it Case 1: $x=t$ and $s\in A$.} The transitivity (on $r, s, u_2$)
implies $r\rightarrow u_2$.  Suppose that $r$ sees $u_3$.  The
transitivity (on $s, r, u_3$) implies $r\rightarrow u_3$ and (on $r,
u_3, u_4$) $r\rightarrow u_4$.  Then $r$ sees $u_5$, for otherwise
$x$-$u_3ru_4$-$u_5$ is a second bull; and $r$ sees $u_6$, for
otherwise $u_2$-$ru_4u_5$-$u_6$ is a second bull; but then
$x$-$u_3u_4r$-$u_6$ is a second bull.  So $r$ misses $u_3$.  Then $r$
sees $u_4$, for otherwise $r$-$u_2xu_3$-$u_4$ is a second bull.  Then
$s$ sees $u_4$, for otherwise $r, s, x, u_3, u_4$ induce a $C_5$.
Then $r$ sees $u_5$, for otherwise $x$-$sru_4$-$u_5$ is a second bull;
and $r$ sees $u_6$, for otherwise $u_3$-$u_4ru_5$-$u_6$ is a second
bull.  But then $x$-$su_4r$-$u_6$ is a second bull.   

{\it Case 2: $x=t$ and $s=u_3$.}  The transitivity (on $r, u_2, u_3, 
u_4$) implies $r\rightarrow u_2$ and $r\rightarrow u_4$.  Then $r$ 
sees $u_5$, for otherwise $x$-$u_3ru_4$-$u_5$ is a second bull; and 
$r$ sees $u_6$, for otherwise $u_2$-$ru_4u_5$-$u_6$ is a second bull. 
But then $x$-$u_3u_4r$-$u_6$ is a second bull. 

{\it Case 3: $x=r$ and $s\in B$.} The transitivity (on $u_3, s, t$)
implies $u_3\rightarrow t$.  If $t$ misses $u_2$, then it sees $u_1$,
for otherwise $u_1$-$u_2xu_3$-$t$ is a second bull; but then $u_1,
u_2, x, s, t$ induce a $C_5$.  So $t$ sees $u_2$.  The transitivity 
(on $s, t, u_2$) implies $u_2\rightarrow t$, and (on $u_1, u_2, t$) 
$u_1\rightarrow t$.   Then $t$ misses $u_4$, for otherwise 
$x$-$u_2u_1t$-$u_4$ is a second bull.  But then $u_1$-$tsu_3$-$u_4$ 
is a second bull.  

{\it Case 4: $x=r$ and $s=u_2$.} The transitivity (on $u_1, u_2, t$)
implies $u_1\rightarrow t$, and similarly we have $u_3\rightarrow t$.
then $t$ misses $u_j$ with $j\in\{4, 5\}$, for otherwise
$x$-$u_2u_1t$-$u_j$ is a second bull; and $t$ misses $u_6$, for 
otherwise $t, u_3, u_4, u_5, u_6$ induce a $C_5$.  But now $t$ is in 
$U_2$, and the fact that $u_2\rightarrow t$ contradicts (\ref{inu2}). 
This completes the proof of the lemma.  $\Box$

\subsection*{Proof of Theorem~\ref{to}}%\footnote{As done in Rio in July.}

The proof of Theorem~\ref{to} goes by induction on the total number of
sensitive vertices in $G$.  We distinguish between two parts, (I) and
(II).

(I) {\it First suppose that $G$ has no sensitive vertex.}
%cyclic bull and no pericyclic bull.} 
By Theorem~\ref{nocycb}, $G$ admits a box partition.
Consider any box $B$.  If $B$ contains any graph $F_j$ with $j=1,
%\ldots, 4$, then, using the auxiliary vertices $a_B, b_B$, we find a
2, 3$, then, using the auxiliary vertices $a_B, b_B$, we find a
spiked $F_j$, which contradicts the fact that $G$ is in ${\cal B}$.
So $B$ contains no bull and no lock.  Gallai~\cite{gal67,mafpre} gave
the list of all minimal forbidden subgraphs for the class of
transitively orientable graphs.  It is a routine matter to check that
every forbidden subgraph in Gallai's list contains either an antihole,
or a bull or a lock.  It follows that every box $B$ induces a subgraph
that admits a transitive orientation $TO(B)$.  Now we make an
orientation of the edges of $G$ by applying the rules below.  In these
rules we use the notation $u\rightarrow v$ to denote the orientation
of an edge $uv$ from $u$ to $v$.  Let us say that an edge $uv$ in a
box $B$ is \emph{sharp} if there is a vertex of $P(B)$ that sees
exactly one of $u,v$, and \emph{dull} otherwise.
\begin{itemize}
\item
Rule $0$: If $uv$ is an edge where $u$ is an odd vertex and $v$ is an
even vertex, then put $u\rightarrow v$.
\item
Rule $S$: If $uv$ is a sharp edge in an odd box $B$, and there is a
vertex of $P(B)$ that sees $u$ and misses $v$, then put $u\rightarrow
v$.  In an even box, put $v\rightarrow u$.
\item
Rule $P3$: If $uv$ is a dull edge in an odd box $B$, and there is a
chordless path $u$-$v$-$w$ in $B$ and a vertex of $P(B)$ that sees $w$
and misses $u, v$, then put $u\rightarrow v$.  In an even box, put
$v\rightarrow u$.
\item
Rule $P4$: If $uv$ is a dull edge in an odd box $B$, and there is a
chordless path $u$-$v$-$w$-$z$ in $B$ and a vertex of $P(B)$ that sees
$z$ and misses $u, v, w$, and $vw$ is dull, then put $v\rightarrow u$.
In an even box, put $u\rightarrow v$.
\item
Rule $Q3$: If $uv$ is a dull edge in an odd box $B$, and there is a
chordless path $u$-$v$-$q$ in $B$ and a vertex of $P(B)$ that sees
$u,v$ and misses $q$, then put $v\rightarrow u$.  In an even box, put
$u\rightarrow v$.
\item
Rule $Q4$: If $uv$ is a dull edge in an odd box $B$, and there is a
chordless path $u$-$v$-$q$-$r$ in $B$ and a vertex of $P(B)$ that sees
$u, v, q$ and misses $r$, then put $u\rightarrow v$.  In an even box,
put $v\rightarrow u$.
\item
Rule $D$: If a dull edge in a box $B$ has not been oriented by the
preceding rules, then orient it according to $TO(B)$.
\end{itemize}
Note that the rules give a symmetric role to odd boxes and even boxes.
Let us prove that these rules produce a transitive orientation of $G$.
\begin{claim}\label{ok}
Every edge of $G$ receives exactly one orientation.
\end{claim}
%
% \subsection*{Proof of claim \ref{ok} in Theorem \ref{to}}
Clearly, Rules $0$, $S$ and $D$ imply that every edge receives at
least one orientation.  Suppose that some edge $uv$ receives the two
opposite orientations $u\rightarrow v$ and $v\rightarrow u$.  By Rule
$0$, edge $uv$ is not between two boxes.  If $uv$ is a sharp edge, the
opposite orientations must both be caused by Rule $S$, so there is a
vertex of $P(B)$ that sees $u$ and misses $v$ and a vertex of $P(B)$
that sees $v$ and misses $u$; but this contradicts
Lemma~\ref{lem:moreitems}.  So $uv$ is a dull edge, say in an even
box.  It cannot be oriented in two opposite ways by Rule $D$, so each
of the two opposite orientations is caused by Rules $P3, P4, Q3, Q4$.
Up to symmetry this yields ten cases, which we analyse now.  In either
case we can consider the auxiliary vertices $a_B, b_B$ for $B$.
Suppose that the two opposite orientations are caused by: \\
- $P3$ and $P3$: So there is a chordless path $u$-$v$-$w$-$x$ with
$w\in B$ and $x\in P(B)$, and there is a chordless path
$v$-$u$-$z$-$y$ with $z\in B$ and $y\in P(B)$.  If $z$ misses $w$,
then $x$ misses $z$, for otherwise $x, z, u, v, w$ induce a $C_5$, and
similarly $y$ misses $w$; but then $x$-$wva_B$-$z$ and $y$-$zua_B$-$w$
are two intersecting bulls.  So $z$ sees $w$.  Then one of $xz, yw$ is
an edge, for otherwise $zw$ is an edge that would be oriented in two
opposite ways by Rule $S$, a contradiction.  Say $x$ sees $z$.  Then
$b_B$ misses $x$, for otherwise $b_B$-$xzw$-$v$ and $b_B$-$xwz$-$u$
are two intersecting bulls.  Then $b_B$-$a_Bvw$-$x$ is a bull.  Then
$b_B$ sees $y$, for otherwise $b_B$-$a_Buz$-$y$ is a second bull.
Then $y$ misses $w$, for otherwise $b_B$-$ywz$-$u$ is a second bull.
Then $y$ sees $x$, for otherwise $y$-$zxw$-$v$ is a second bull.  But
then $b_B$-$yxz$-$u$ is a second bull, a contradiction. \\
- $P3$ and $P4$: So there is a chordless path $u$-$v$-$w$-$x$ with
$w\in B, x\in P(B)$, and there is a chordless path $u$-$v$-$s$-$t$-$y$
with $s, t\in B, y\in P(B)$, and $vs$ is dull.  So $x$ misses $s$.
Then $u$-$a_Bst$-$y$ is a bull.  Then $w$ misses $s$, for otherwise
$u$-$vsw$-$x$ is a second bull.  But then $vs$ is oriented in two
opposite ways by Rules $P3$ and $P3$ (because of $x$-$w$-$v$-$s$ and
$y$-$t$-$s$-$v$), a contradiction.  \\
- $P3$ and $Q3$: So there is a path $u$-$v$-$w$-$x$ with $w\in B, x\in
P(B)$ and a path $u$-$v$-$q$ with $q\in B$ and a vertex $y$ that sees
$u,v$ and misses $q$.  Note that either $b_B$-$a_Bqv$-$y$ or
$b_B$-$yuv$-$q$ is one bull.  Then $y$ sees $w$, for otherwise $vw$ is
a sharp edge oriented both ways by Rule $S$ (because of $x,y$), which
contradicts a fact already proved.  Then $x$ misses $q$, for otherwise
$vq$ is a sharp edge oriented both ways (because of $x,y$).  Then $x$
misses $y$, for otherwise $x$-$yuv$-$q$ is a second bull.  Then $w$
sees $q$, for otherwise $x$-$wyv$-$q$ is a second bull.  But then
$x$-$wqv$-$u$ is a second bull, a contradiction. \\
- $P3$ and $Q4$: So there is a chordless path $u$-$v$-$w$-$x$ with
$w\in B, x\in P(B)$, and there is a chordless path $v$-$u$-$q$-$r$
with $q,r\in B$ and a vertex $y\in P(B)$ that sees $v, u, q$ and
misses $r$, and $uq$ is dull.  So $x$ misses $q$.  Note that either
$b_B$-$a_Brq$-$y$ or $b_B$-$yuq$-$r$ is one bull that contains $q$.
Then $x$ misses $r$, for otherwise $rq$ is a sharp edge oriented both
ways (because of $x,y$).  Then $y$ sees $w$, for otherwise $vw$ is a
sharp edge oriented both ways (because of $x,y$).  Then $w$ misses
$q$, for otherwise $uq$ is oriented both ways by $P3$ (because of
$x$-$w$-$q$-$u$) and $Q3$ (because $u$-$q$-$r$ and $y$), which
contradicts a fact already proved.  Then $x$ sees $y$, for otherwise
$x$-$wvy$-$q$ is a second bull.  But then $x$-$yuq$-$r$ is a second
bull, a contradiction.  \\
- The remaining six cases ($P4$ and $P4$; $P4$ and $Q3$; $P4$ and
$Q4$; $Q3$ and $Q3$; $Q3$ and $Q4$; $Q4$ and $Q4$) can all be treated
as follows.  When $u\rightarrow v$ is given by Rule $P4$, there is a
chordless path $u$-$v$-$w$-$z$-$x$ with $w,z\in B, x\in P(B)$, and
then $u$-$a_Bwz$-$x$ is a bull.  When $u\rightarrow v$ is given by
Rule $Q3$, there is a path $u$-$v$-$q$ in $B$ and some $y\in P(B)$
that sees $u,v$ and misses $q$; then either $b_B$-$a_Bqv$-$y$ or
$b_B$-$yuv$-$q$ is a bull.  When $u\rightarrow v$ is given by Rule
$Q4$, there is a path $v$-$u$-$q$-$r$ in $B$ and a vertex $y\in P(B)$
that sees $v, u, q$ and misses $r$; then either $b_B$-$a_Brq$-$y$ or
$b_B$-$yuq$-$r$ is a bull.  And so when $v\rightarrow u$ is given by
Rules $P4, Q3, Q4$, there is a similar bull.  It is a routine matter
to check that in each of the six cases, the two bulls produced by the
two rules are distinct and intersect, a contradiction.  Therefore
(\ref{ok}) holds.
\begin{claim}\label{trans}
The orientation produced by the rules is transitive.
\end{claim}
%
% \subsection*{Proof of claim \ref{trans} in Theorem \ref{to}}
Consider any chordless path $u$-$v$-$w$ in $G$.  Assume that $u,v,w$
are not all in the same box.  Then, up to symmetry, one of them is odd
and the other two are even (or vice-versa).  If $v$ is the odd one, we
have $v\rightarrow u$ and $v\rightarrow w$ by Rule $0$, so $u$-$v$-$w$
is oriented transitively.  If $u$ is the odd one, we have
$u\rightarrow v$ by Rule $0$ and $w\rightarrow v$ by Rule $S$, so
$u$-$v$-$w$ is oriented transitively.  Now we may assume up to
symmetry that $u,v,w$ are all in one odd box $B$.  If both $uv, vw$
are oriented by Rule $D$, then $u$-$v$-$w$ is oriented transitively
since $TO(B)$ is a transitive orientation.  So we may assume that at
least one of $uv, vw$, say $uv$, is oriented by one of Rules $S, P3,
P4, Q3, Q4$.  Suppose by contradiction that the rules produce
$u\rightarrow v$ and $v\rightarrow w$.  In either case we can consider
the auxiliary vertices $a_B, b_B$ for $B$.  Let us analyze all the
cases.  \\
- Suppose that the two orientations $u\rightarrow v$ and $v\rightarrow
w$ are both caused by $S$.  So there is a vertex $x\in P(B)$ that sees
$u$ and misses $v$, and there is a vertex $y\in P(B)$ that sees $v$
and misses $w$.  Then $y$ sees $u$, for otherwise $uv$ is oriented in
two opposite ways by Rule $S$ (because of $x,y$), which contradicts
(\ref{ok}).  Then $x$ misses $w$, for otherwise $vw$ is oriented both
ways by $S$ (because of $x,y$).  Note that either $b_B$-$awv$-$y$ (if
$b_B$ misses $y$) or $b_B$-$yuv$-$w$ (if $b_B$ sees $y$) is one bull.
So $x$ sees $y$, for otherwise $x$-$uyv$-$w$ is a second bull.  But
now, vertices $u, v, w, x, y$ contradict Property (vii) for $B$.  \\
So we may now assume, up to symmetry, that $vw$ is dull.  \\
- Suppose that $u\rightarrow v$ is caused by $S$.  So there is a
vertex $x\in P(B)$ that sees $u$ and misses $v$.  Then $x$ misses $w$
since $vw$ is dull.  But then $w\rightarrow v$ is given by Rule $P3$
(because of $x$-$u$-$v$-$w$), which contradicts (\ref{ok}).  \\
- Suppose that $u\rightarrow v$ is caused by $P3$.  So there is a
chordless path $u$-$v$-$z$-$x$ with $z\in B, x\in P(B)$.  Then $x$
misses $w$ since $vw$ is dull.  If $z$ sees $w$, then $x$-$zwv$-$u$
and $x$-$zwa_B$-$u$ are two intersecting bulls, a contradiction.  So
$z$ misses $w$.  Then $w\rightarrow v$ is given by Rule $P3$ (because
of $w$-$v$-$z$-$x$), which contradicts (\ref{ok}).  \\
- Suppose that $u\rightarrow v$ is caused by $P4$.  So there is a
chordless path $v$-$u$-$z$-$p$-$x$ with $z, p\in B, x\in P(B)$, and
$zu$ is dull.  So $x$ misses $w$, since $vw$ is dull.  Note that
$x$-$pza_B$-$v$ is one bull.  If $p$ sees $w$, then $z$ sees $w$, for
otherwise $p, z, u, v, w$ induce a $C_5$; but then $x$-$pzw$-$v$ is a
second bull.  So $p$ misses $w$.  Then $z$ sees $w$, for otherwise
$x$-$pza_B$-$w$ is a second bull.  But then $w\rightarrow v$ is given
by Rule $P4$ (because of $x$-$p$-$z$-$w$-$v$), which contradicts
(\ref{ok}).  \\
- Suppose that $u\rightarrow v$ is caused by $Q3$.  So there is a
chordless path $v$-$u$-$q$ in $B$ and a vertex $y\in P(B)$ that sees
$v,u$ and misses $q$.  So $y$ misses $w$, since $vw$ is dull.  If $q$
sees $w$, then $w\rightarrow v$ is given by Rule $Q3$ (because of
$v$-$w$-$q$ and $y$), a contradiction.  If $q$ misses $w$, then
$w\rightarrow v$ is given by Rule $Q4$ (because of $w$-$v$-$u$-$q$ and
$y$), a contradiction.  \\
- Finally suppose that $u\rightarrow v$ is caused by $Q4$.  So there
is a chordless path $u$-$v$-$q$-$r$ in $B$ and a vertex $y\in P(B)$
that sees $u, v, q$ and misses $r$.  Then $y$ sees $w$ since $vw$ is
dull.  If $w$ sees $r$, then $w\rightarrow v$ is given by Rule $Q3$
(because of $v$-$w$-$r$ and $y$), a contradiction.  So $w$ misses $r$.
If $w$ misses $q$, then $w\rightarrow v$ is given by Rule $Q4$
(because of $w$-$v$-$q$-$r$ and $y$), a contradiction.  So $w$ sees
$q$.  But then $r$-$qwv$-$u$ and $r$-$qwy$-$u$ are two intersecting
bulls.  Therefore (\ref{trans}) holds.

A classical theorem of Ghouila-Houri~\cite{gho64} states that if a
graph admits a transitive orientation then it admits a transitive and
acyclic orientation.  So (\ref{trans}) suffices to prove our theorem
and the proof of part (I) is complete.  Actually, it is not hard to
prove that the orientation produced by the above rules has no circuit,
but we omit this proof.

(II) {\it Now $G$ has a sensitive vertex.} Let $x$ be any sensitive
vertex of $G$.  By the definition of a sensitive vertex, $G-x$
contains a hole of length at least six.  Thus $G-x$ is in class ${\cal
B}$ and contains a hole.  By the induction hypothesis, $G-x$ has a
transitive orientation; and by Lemma~\ref{lem:sensto}, $G$ has a
transitive orientation.  This completes the proof of Theorem~\ref{to}.
$\Box$

\section{A colouring algorithm}
\label{s:co}

We conclude the paper with a discussion about how Theorem~\ref{main1}
indeed yields a polynomial-time algorithm that colours the vertices of
any bull-reducible Berge graph containing no antihole.

% Let us discuss the complexity of our algorithm.  
We are given a bull-reducible Berge graph $G$ with no antihole, 
with $n$ vertices and
$m$ edges.

In the preliminary step, we will use the algorithm of Spinrad
{\cite{spi-homo}}, which finds all maximal homogeneous sets of a
graph.  The complexity of Spinrad's algorithm is O$(mf(n,m))$, where
$f(n,m)$ is the reverse of the Ackerman function.  Remark that the
maximal homogeneous sets are pairwise disjoint.  For each such
homogeneous set $H$, we can apply recursively our algorithm on $H$ and
find a coloring of $H$ with $\omega(H)$ colors.  Then we replace in
$G$ the vertices of $H$ by a clique $Q(H)$ of size $\omega(H)$, and do
this for each maximal homogeneous set.  Trivially the resulting graph
is isomorphic to a subgraph of the original graph.  At the end, it is
easy to get a coloring of the original graph from a coloring of the
new graph simply by merging the colors used on $Q(H)$ with the colors
used in $H$.

In the second step, we determine whether the graph is weakly
triangulated using the following ``naive'' method.  For each triple
$abc$ forming a $P_3$ we test if this $P_3$ extends to a hole in the
graph.  Clearly, it suffices to check whether $a$ and $c$ are in the
same component of the subgraph obtained from $G$ by removing the
vertices in $N(a)\cap N(c)$ and the vertices in $N(b)-\{a, c\}$.
Using a shortest path algorithm we will find a shortest hole
containing $a,b,c$, if any.  Globally, we will either find that $G$ is
weakly triangulated or determine a shortest even hole in $G$.

If $G$ is weakly triangulated, we refer to the algorithm in~\cite{hhm89}.

In the remaining case, if $G$ has no sensitive vertex then $G$ admits
a box partition which can be determined by breadth-first search.  Now,
we may apply the rules to the box partition and obtain a transitive
orientation for $G$.  Then we apply the greedy method on the
transitive orientation.  Else, if $G$ has a sensitive vertex $x$, then
Lemma~\ref{lem:sensto} extends a transitive orientation from $G-x$ to
$G$.

The overall complexity is O$(n^4m)$.

\subsection*{The weighted case}

Let us remark that this algorithm can be adapted to solve the weighted
version of the coloring problem.  Given a graph with a weight function
$w$ on its vertices, a {\it weighted coloring\/} is a family of stable
sets $S_1,\ldots, S_q$ with weights $W(S_i)$ such that:
\begin{equation}
w(x)\le \sum_{S_i\ni x} W(S_i)\label{eq:weicol}
\end{equation}
holds for every vertex $x$.  The goal is then to find a weighted
coloring whose total weight $W(S_1)+\cdots+W(S_q)$ is minimal.  With a
our algorithm we can solve the minimum-weight coloring problem as
follows.
\newline
If the graph has an incomplete homogeneous set $H$, we can recursively
apply the algorithm on $H$ and find a minimum-weight coloring for $H$.
This consists of a family of stable sets $T_1,\ldots, T_p$ with
weights $W(T_1),\ldots, W(T_p)$.  We then substitute $H$ by a clique
$Q(H)$ of cardinality $p$, resulting in a new graph $G'$.  The $i$-th
vertex in $Q(H)$ receives weight $W(T_i)$, while the vertices of
$G'-Q$ keep the same weight as in the original graph.  Classical
polyhedral considerations (see \cite{grolovsch84}) imply that in a
weighted perfect graph there exists a minimum-weight coloring that
satisfies (\ref{eq:weicol}) with equality for every vertex.  So, at
the end, we can obtain a minimum-weight coloring of the original graph
by combining a minimum-weight of the new graph with $T_1,\ldots, T_p$.
\newline
When the graph is weakly triangulated, we can apply the minimum-weight
coloring algorithm from {\cite{hhm89}}.
\newline
When the graph is transitively orientable, we can apply the minimum-weight
coloring algorithm from {\cite{hoa90cs}}, whose complexity is O$(nm)$.

The overall complexity is again O$(n^4m)$.

To find a maximum-weight clique is straightforwardly similar, because
the algorithms respectively from {\cite{hhm89}} and {\cite{hoa90cs}}
can also be required to produce a maximum weighted clique in
respectively a weakly triangulated graph and a transitively orientable
graph.

\clearpage

\end{document}